\documentclass{article}

\usepackage[final, nonatbib]{nips_2017} 

\usepackage{amsmath,amssymb,amsthm,amsfonts}
\usepackage[utf8]{inputenc} 
\usepackage[T1]{fontenc}    
\usepackage{color}
\usepackage{algpseudocode}
\usepackage{algorithm,algorithmicx}
\usepackage{hyperref,url}       
\usepackage{booktabs}       
\usepackage{nicefrac}       
\usepackage{microtype}      
\usepackage[numbers,sort&compress]{natbib}
\usepackage{caption}

\definecolor{medblue}{rgb}{0,0,.75}
\algrenewcommand\alglinenumber[1]{\sf\tiny\color{medblue}{#1}\quad}
\algrenewcommand\algorithmicrequire{\textbf{Input:}}
\algrenewcommand\algorithmicensure{\textbf{Output:}}
\newcommand{\R}{\mathbb{R}}
\newcommand{\Rext}{\mathbb{R}\cup\{+\infty\}}
\newcommand{\set}[1]{\left\{#1\right\}}

\newcommand{\norm}[1]{\left\Vert#1\right\Vert}
\newcommand{\dom}[1]{\mathrm{dom}\left(#1\right)}
\newcommand{\prox}[2]{\mathrm{prox}_{#1}\left(#2\right)}
\newcommand{\argmin}{\mathrm{arg}\min}
\newcommand{\Eproof}{\hfill $\square$}
\newcommand{\Fc}{\mathcal{F}}

\newcommand{\Xc}{\mathcal{X}}

\newcommand{\Id}{\mathbb{I}}

\newcommand{\iprods}[1]{\langle #1\rangle}

\newcommand{\prob}[1]{\mathbb{P}\left\{#1\right\}}
\newcommand{\expect}[2]{\mathbb{E}_{#1}\left[#2\right]}

\newtheorem{theorem}{Theorem}[section]

\newtheorem{corollary}[theorem]{Corollary}
\newtheorem{lemma}[theorem]{Lemma}
\newtheorem{proposition}[theorem]{Proposition}
\newtheorem{remark}[theorem]{Remark}
\newtheorem{assumption}{Assumption}

\usepackage{graphicx} 
\usepackage{subfigure}

\title{Smooth Primal-Dual Coordinate Descent Algorithms for Nonsmooth Convex Optimization}

\author{
Ahmet Alacaoglu\textsuperscript{1} 
\And
Quoc Tran-Dinh\textsuperscript{2} 
\And
Olivier Fercoq\textsuperscript{3} 
\And
Volkan Cevher\textsuperscript{1} 
\AND
\textsuperscript{1}{\normalfont Laboratory for Information and Inference Systems (LIONS), EPFL, Lausanne, Switzerland} \\
  \texttt{\{ahmet.alacaoglu, volkan.cevher\}@epfl.ch} \\
  \textsuperscript{2} {\normalfont Department of Statistics and Operations Research, UNC-Chapel Hill, NC, USA} \\
\texttt{quoctd@email.unc.edu} \\
  \textsuperscript{3} {\normalfont LTCI, Télécom ParisTech, Université Paris-Saclay, Paris, France} \\
\texttt{olivier.fercoq@telecom-paristech.fr}
}

\begin{document}

\maketitle
\begin{abstract}
\vspace{-2ex}
We propose a new randomized coordinate descent method for a convex optimization template 
with broad applications. Our analysis relies on a novel combination of four
ideas applied to the primal-dual gap function: smoothing, acceleration, homotopy, and coordinate descent with non-uniform sampling. As a result, our method features the first convergence rate guarantees among the coordinate descent methods, that are the best-known under a variety of common structure assumptions on the template. We provide numerical evidence to support the theoretical results with a comparison to state-of-the-art algorithms. 
\end{abstract}

\vspace{-2ex}
\section{Introduction}\label{sec:intro}
\vspace{-2ex}
We develop randomized coordinate descent methods to solve the following composite convex problem:
\begin{align}\label{eq:cvx_prob}
F^{\star} = \min_{x \in\R^p} \set{F(x) = f(x) + g(x) + h(Ax) }, 
\end{align}
where $f : \R^p \to \R$, $g : \R^p \to \Rext$, and $h : \R^m\to\Rext$ are proper, closed and convex functions, $A \in \R^{m \times p}$ is a given matrix.  
The optimization template~\eqref{eq:cvx_prob} covers many important applications including support vector machines, sparse model selection, logistic regression, etc. It is also convenient to formulate generic constrained convex problems by choosing an appropriate  $h$. 

\vspace{-0.5ex}
Within convex optimization, coordinate descent methods have recently become increasingly popular in the literature \cite{Nesterov2012,richtarik2014iteration,richtarik2016parallel,fercoq2015accelerated,shalev2013stochastic,necoara2016parallel}.
These methods are particularly well-suited to solve huge-scale problems arising from machine learning applications where matrix-vector operations are prohibitive \cite{Nesterov2012}.  

\vspace{-0.5ex}
To our knowledge, there is no coordinate descent method for the general three-composite form  \eqref{eq:cvx_prob} within our structure assumptions  studied here that has rigorous convergence guarantees. 
Our paper specifically fills this gap. 
For such a theoretical development, coordinate descent algorithms require specific assumptions on the convex optimization problems \cite{Nesterov2012,fercoq2015accelerated,necoara2016parallel}. 
As a result, to rigorously handle the three-composite case, we assume that ($i$) $f$ is smooth, ($ii$) $g$ is non-smooth but decomposable (each component has an ``efficiently computable'' proximal operator), and ($iii$) $h$ is non-smooth.

\vspace{-2ex}
\paragraph{Our approach:} In a nutshell, we generalize \cite{fercoq2015accelerated,qu2016coordinate} to the three composite case \eqref{eq:cvx_prob}. For this purpose, we combine several classical and contemporary ideas: We exploit the  smoothing technique in \cite{Nesterov2005c}, the efficient implementation technique in \cite{lee2013efficient,fercoq2015accelerated}, the homotopy strategy in \cite{tran2015smooth}, and the nonuniform coordinate selection rule in \cite{qu2016coordinate} in our algorithm, to achieve the best known complexity estimate for the template. 

Surprisingly, the combination of these ideas is achieved in a very natural and elementary primal-dual gap-based framework. However, the extension is indeed not trivial since it requires to deal with the composition of a non-smooth function $h$ and a linear operator $A$. 

While our work has connections to the methods developed in \cite{qu2016coordinate,fercoq2013smooth,fercoq2015coordinate}, it is rather distinct. 
First, we consider a more general problem \eqref{eq:cvx_prob} than the one in \cite{fercoq2015accelerated, qu2016coordinate,fercoq2013smooth}. 
Second, our method relies on Nesterov's accelerated scheme rather than a primal-dual method as in \cite{fercoq2015coordinate}. 
Moreover, we obtain the first rigorous convergence rate guarantees as opposed to \cite{fercoq2015coordinate}. In addition, we allow using any sampling distribution for choosing the coordinates.

\textbf{Our contributions:} 
We propose a new smooth primal-dual randomized coordinate descent method for solving \eqref{eq:cvx_prob} where $f$ is smooth, $g$ is nonsmooth, separable and has a block-wise proximal operator, and $h$ is a general nonsmooth function. 
Under such a structure, we show that our algorithm achieves the best known $\mathcal{O}(n/k)$ convergence rate, where $k$ is the iteration count and to our knowledge, this is the first time that this convergence rate is proven for a coordinate descent algorithm. 

We  instantiate our algorithm to solve special cases of \eqref{eq:cvx_prob} including the case $g=0$ and constrained problems.
We analyze the convergence rate guarantees of these variants individually and discuss the choices of sampling distributions.

Exploiting the strategy in \cite{lee2013efficient, fercoq2015accelerated}, our algorithm can be implemented in parallel by breaking up the full vector updates. We also provide a restart strategy to enhance practical performance. 

\vspace{-2ex}
\paragraph{Paper organization:}
We review some preliminary results in Section~\ref{sec:preliminary}. 
The main contribution of this paper is in Section~\ref{sec:main_alg} with the main algorithm and its convergence guarantee.
We also present special cases of the proposed algorithm. 
Section~\ref{sec:numerical_experiments} provides numerical evidence to illustrate the performance of our algorithms in comparison to existing methods.
The proofs are deferred to the supplementary document.

\vspace{-2ex}
\section{Preliminaries}\label{sec:preliminary}
\vspace{-2.5ex}
\paragraph{Notation:}
Let $[n] := \set{1, 2, \cdots, n}$ be the set of $n$ positive integer indices. 
Let us decompose the variable vector $x$ into $n$-blocks denoted by $x_i$ as $x = [x_1; x_2; \cdots; x_n]$ such that each block $x_i$ has the size $p_i \geq 1$ with $\sum_{i=1}^np_i = p$.
We also decompose the identity matrix $\Id_p$ of $\mathbb{R}^p$ into $n$ block as $\Id_p = [U_1, U_2, \cdots, U_n]$, where $U_i \in \mathbb{R}^{p \times p_i}$ has $p_i$ unit vectors.
In this case, any vector $x\in\R^p$ can be written as $x = \sum_{i=1}^n U_i x_i$, and each block becomes $x_i = U_i^\top x$ for $i \in [n]$. We define the partial gradients as $\nabla_i f(x) = U_i^\top \nabla f(x)$ for $i \in [n]$. 
For a convex function $f$, we use $\dom{f}$ to denote its domain, $f^{\ast}(x) := \sup_{u}\set{u^{\top}x - f(u)}$ to denote its Fenchel conjugate, and $\mathrm{prox}_f(x) := \argmin_u\set{ f(u) + (1/2)\Vert u -x\Vert^2 }$ to denote its proximal operator.
For a convex set $\Xc$, $\delta_{\Xc}(\cdot)$ denotes its indicator function.
\vspace{-0.75ex}
We also need the following weighted norms:
\vspace{-0.5ex}
\begin{equation}\label{eq:norms}
{\!\!\!\!\!}\begin{array}{llllll}
\Vert x_i \Vert^2_{(i)} &{\!\!\!\!\!\!\!}= \iprods{H_ix_i,x_i}, &{\!\!\!\!}\text{}~~~~(\Vert y_i \Vert_{(i)}^{\ast})^2  &{\!\!\!\!\!\!\!}= \iprods{H_i^{-1}y_i,y_i},  \vspace{1ex}\\
\Vert x\Vert_{[\alpha]}^2  &{\!\!\!\!\!\!\!}= \sum_{i=1}^nL_i^{\alpha}\Vert x_i \Vert^2_{(i)}, &{\!\!\!\!}\text{}~~~~(\Vert y \Vert_{[\alpha]}^{\ast})^2 &{\!\!\!\!\!\!\!}= \sum_{i=1}^nL_i^{-\alpha}(\Vert y_i\Vert_{(i)}^\ast)^2.
\end{array}{\!\!\!}
\vspace{-0.5ex}
\end{equation}

Here, $H_i\in\R^{p_i\times p_i}$ is a symmetric positive definite matrix, and $L_i\in (0, \infty)$ for $i\in [n]$ and $\alpha > 0$. In addition, we use $\Vert \cdot \Vert$ to denote $\Vert \cdot \Vert _2$.

\vspace{-2ex}
\paragraph{Formal assumptions on the template:}
We require the following assumptions to tackle \eqref{eq:cvx_prob}: 
\begin{assumption}\label{as:fun_assumption}
The functions $f$, $g$ and $h$ are all proper, closed and convex. Moreover, they satisfy 
\begin{itemize}
\vspace{-0.75ex}
\item[(a)] The partial derivative $\nabla_i{f}(\cdot)$ of $f$ is Lipschitz continuous with the Lipschitz constant $\hat{L}_i \in [0, +\infty)$, i.e., 
$\Vert \nabla_if(x + U_id_i) - \nabla_if(x)\Vert_{(i)}^{\ast} \leq \hat{L}_i\Vert d_i\Vert_{(i)}$ for all $x\in\R^p, d_i\in\R^{p_i}$.
\vspace{-1ex}
\item[(b)] The function $g$ is separable, which has the following form $g(x) = \sum_{i=1}^n g_i(x_{i})$.
\vspace{-1ex}
\item[(c)] One of the following assumptions for $h$ holds  for Subsections~\ref{subsec:case1} and~\ref{sec:alg_constrained}, respectively:
\begin{itemize} 
\vspace{-1.5ex}
\item[i.] \label{assumption:h_lips} $h$ is Lipschitz continuous which is equivalent to the boundedness of $\dom{h^{\ast}}$.
\item[ii.] \label{assumption:h_ind} $h$ is the indicator function for an equality constraint, i.e., $h(Ax) := \delta_{\{c\}}(Ax)$.
\end{itemize}
\vspace{-1ex}
\end{itemize}
\end{assumption}
Now, we briefly describe the main techniques used in this paper.
\vspace{-2ex}
\paragraph{Acceleration:} 
Acceleration techniques in convex optimization date back to the seminal work of Nesterov in \cite{nesterov1983method}, and is one of standard techniques in convex optimization.
We exploit such a scheme to achieve the best known $\mathcal{O}(1/k)$ rate for the nonsmooth template \eqref{eq:cvx_prob}.
\vspace{-2ex}
\paragraph{Nonuniform distribution:}
We assume that $\xi$ is a random index on $[n]$ associated with a probability distribution $q = (q_1,\cdots, q_n)^\top$ such that 
\begin{equation}\label{eq:prob_i}
\prob{\xi = i} = q_i > 0,~~i\in [n],~~~\text{and}~~\sum_{i=1}^nq_i = 1.
\end{equation} 
When $q_i = \frac{1}{n}$ for all $i\in[n]$, we obtain the uniform distribution.
Let $i_0, i_1, \cdots, i_k$ be i.i.d. realizations of the random index $\xi$ after $k$ iteration. We define $\Fc_{k+1} = \sigma(i_0, i_1,\cdots, i_k)$ as the $\sigma$-field generated by these realizations.

\vspace{-2ex}
\paragraph{Smoothing techniques:}
We can write the convex function $h(u) = \sup_{y}\set{\iprods{u, y} - h^{\ast}(y)}$ using its Fenchel conjugate $h^{\ast}$. 
Since  $h$ in \eqref{eq:cvx_prob} is convex but possibly nonsmooth, we smooth $h$ as
\begin{equation}\label{eq:smooth_h}
h_{\beta}(u) := \max_{y\in\R^m}\set{ \iprods{u, y} - h^{\ast}(y) - \tfrac{\beta}{2}\Vert y - \dot{y}\Vert^2},
\end{equation}
where $\dot{y}\in\R^m$ is given and $\beta > 0$ is the smoothness parameter.
Moreover, the quadratic function $b(y, \dot{y}) = \frac{1}{2}\Vert y - \dot{y}\Vert^2$ is defined based on a given norm in $\R^m$.
Let us denote by $y^{\ast}_{\beta}(u)$, the unique solution of this concave maximization problem in \eqref{eq:smooth_h}, i.e.:
\begin{equation}\label{eq:y_sol}
y^{\ast}_{\beta}(u) := \mathrm{arg}\max_{y\in\R^m}\set{ \iprods{u, y} - h^{\ast}(y) - \tfrac{\beta}{2}\Vert y - \dot{y}\Vert^2} = \prox{\beta^{-1}h^{\ast}}{\dot{y} + \beta^{-1}u},
\end{equation}
where $\text{prox}_{h^\ast}$ is the proximal operator of $h^\ast$. If we assume that $h$ is Lipschitz continuous, or equivalently that $\dom{h^{\ast}}$ is bounded, then it holds that
\begin{equation}\label{eq:bound_est}
h_{\beta}(u) \leq h(u) \leq h_{\beta}(u) + \tfrac{\beta D_{h^{\ast}}^2}{2},~~~~\text{where}~~D_{h^{\ast}} := \max_{y\in\dom{h^{\ast}}} \Vert y - \dot{y}\Vert < +\infty.
\end{equation}
Let us define a new smoothed function $\psi_{\beta}(x) := f(x) + h_{\beta}(Ax)$. Then, $\psi_{\beta}$ is differentiable, and its block partial gradient 
\begin{equation}\label{eq:H_function}
\nabla_i{\psi_{\beta}}(x) = \nabla_if(x) + A_i^\top y^{\ast}_{\beta}(Ax) 
\end{equation}
is also Lipschitz continuous with the Lipschitz constant $L_i(\beta) := \hat{L}_i + \frac{\Vert A_i\Vert^2}{\beta}$, where $\hat{L}_i$ is given in Assumption~\ref{as:fun_assumption}, and $A_i \in\R^{m\times p_i}$ is the $i$-th block of $A$.
\vspace{-2ex}
\paragraph{Homotopy:} In smoothing-based methods, the choice of the smoothness parameter is critical. This choice may require the knowledge of the desired accuracy, number of maximum iterations or the diameters of the primal and/or dual domains as in \cite{Nesterov2005c}. 
In order to make this choice flexible and our method applicable to the constrained problems, we employ a homotopy strategy developed in \cite{tran2015smooth} for deterministic algorithms, to gradually update the smoothness parameter while making sure that it converges to $0$.

\vspace{-2ex}
\section{Smooth primal-dual randomized coordinate descent}\label{sec:main_alg}
\vspace{-2ex}
In this section, we develop a smoothing primal-dual method to solve \eqref{eq:cvx_prob}. 
Or approach is to combine the four key techniques mentioned above: smoothing, acceleration, homotopy, and randomized coordinate descent.
Similar to \cite{qu2016coordinate} we allow to use arbitrary nonuniform distribution, which may allow to design a good distribution that captures the underlying structure of specific problems.

\vspace{-2ex}
\subsection{The algorithm}
\vspace{-2ex}
Algorithm \ref{alg:A1} below smooths, accelerates, and randomizes the coordinate descent method. 
\vspace{-1.5ex}
\begin{algorithm}[H]
\caption{SMooth,  Accelerate, Randomize The Coordinate Descent (SMART-CD)}\label{alg:A1}
\begin{algorithmic}[1]
    \Require{Choose  $\beta_1 > 0$ and $\alpha \in [0, 1]$ as two input parameters. Choose $x^0\in\R^p$}.
    \State Set $B_i^0 := \hat{L}_i + \frac{\Vert A_i\Vert^2}{\beta_1}$ for $i\in [n]$. 
    Compute $S_{\alpha} := \sum_{i=1}^n(B_i^0)^{\alpha}$ and $q_i := \frac{(B_i^0)^{\alpha}}{S_{\alpha}}$ for all $i\in [n]$.
    \State Set $\tau_0 := \min\set{q_i\mid 1\leq i \leq n} \in (0, 1]$ for $i\in [n]$. Set $\bar{x}^0 = \tilde{x}^0 := x^0$.
     \vspace{0.5pc}
    \For{$k \gets 0,1,\cdots, k_{\max}$}
        \State\label{step:x_hat} 
              Update $\hat{x}^k := (1-\tau_k)\bar{x}^k + \tau_k \tilde{x}^k$ and compute $\hat{u}^k := A\hat{x}^k$.
         \State\label{step:y_hat}  
              Compute the dual step $y^{\ast}_k := y^{\ast}_{\beta_{k+1}}(\hat{u}^k) = \prox{\beta_{k+1}^{-1}h^{\ast}}{\dot{y} + \beta_{k+1}^{-1}\hat{u}^k}.$
         \State Select a block coordinate $i_{k} \in [n]$ according to the probability distribution $q$.
	\State\label{step:x_tilde}
	      Set $\tilde{x}^{k+1} := \tilde{x}^k$,  and compute the primal $i_k$-block coordinate:
	      \begin{equation*}
	      \tilde{x}^{k+1}_{i_{k}} := \mathrm{arg}\!\!\!\!\!\!\min_{x_{i_k}\in\R^{p_{i_k}}}\Big\{\iprods{\nabla_{i_k}f(\hat{x}^k) + A_{i_k}^\top y^{\ast}_k, x_{i_k} - \hat{x}^k_{i_k}} + g_{i_k}(x_{i_k}) 
	                                           + \frac{\tau_k B^{k} _{i_{k}}}{2 \tau_0} \Vert x_{i_k} - \tilde{x}^k_{i_{k}} \Vert_{(i_k)}^2 \Big\}.
	     \end{equation*}
           \State\label{step:x_bar}
               Update $\bar{x}^{k+1} := \hat{x}^k + \frac{\tau_k}{\tau_0} (\tilde{x}^{k+1} - \tilde{x}^{k})$.
	   \State\label{step:update_tau} Compute $\tau_{k+1} \in (0, 1)$ as the unique positive root of $\tau^3 + \tau^2 + \tau_k^2 \tau - \tau_k ^2=0$.
	   \State\label{step:update_beta} Update $\beta_{k+2} := \frac{\beta_{k+1}}{1+\tau_{k+1}}$ and $B^{k+1} _i := \hat{L}_i + \frac{\Vert A_i \Vert ^2}{\beta_{k+2}}$ for $i\in [n]$.
     \EndFor
\end{algorithmic}
\end{algorithm}
\vspace{-2ex}
\noindent
From the update $\bar{x}^{k} := \hat{x}^{k-1} + \frac{\tau_{k-1}}{\tau_0} (\tilde{x}^{k} - \tilde{x}^{k-1})$ and $\hat{x}^k := (1-\tau_k)\bar{x}^k + \tau_k \tilde{x}^k$, it directly follows that $\hat{x}^k := (1-\tau_k)\big(\hat{x}^{k-1} + \frac{\tau_{k-1}}{\tau_0} (\tilde{x}^{k} - \tilde{x}^{k-1})\big) + \tau_k \tilde{x}^k$. Therefore, it is possible to implement the algorithm without forming $\bar{x}^{k}$.

\vspace{-2ex}
\subsection{Efficient implementation}
\vspace{-2ex}
While the basic variant in Algorithm~\ref{alg:A1} requires full vector updates at each iteration, we exploit the idea in \cite{lee2013efficient, fercoq2015accelerated} and show that we can partially update these vectors in a more efficient manner.
\vspace{-2ex}
\begin{algorithm}[H]
\caption{Efficient SMART-CD}\label{alg:A_eff}
\begin{algorithmic}[1]
    \Require{Choose a parameter $\beta_1 > 0$ and $\alpha \in [0, 1]$ as two input parameters. Choose $x^0\in\R^p$.}
    \State Set $B_i^0 := \hat{L}_i + \frac{\Vert A_i\Vert^2}{\beta_1}$ for $i\in [n]$. 
    Compute $S_{\alpha} := \sum_{i=1}^n(B_i^0)^{\alpha}$ and $q_i := \frac{(B_i^0)^{\alpha}}{S_{\alpha}}$ for all $i\in [n]$.
    \State Set $\tau_0 := \min\set{q_i\mid 1\leq i \leq n} \in (0, 1]$ for $i\in [n]$ and $c_0 = (1 - \tau_0)$. Set $u^0=\tilde{z}^0 := x^0$.
     \vspace{0.5pc}
    \For{$k \gets 0,1,\cdots, k_{\max}$}
         \State\label{step:y_star}  
              Compute the dual step $y^{\ast}_{\beta_{k+1}}(c_k A u^k + A \tilde{z}^k) := \prox{\beta_{k+1}^{-1}h^{\ast}}{\dot{y} + \beta_{k+1}^{-1} (c_k A u^k + A\tilde{z}^k) }.$
         \State Select a block coordinate $i_{k} \in [n]$ according to the probability distribution $q$.
	\State\label{step:t} Let $\nabla_i^k := \nabla_{i_k}f(c_k u^k + \tilde{z}^{k}) + A_{i_k}^\top y^{\ast}_{\beta_{k+1}}(c_k A u^k + A \tilde{z}^k)$. Compute
	      \begin{equation*}
	      t^{k+1}_{i_{k}} := \mathrm{arg}\!\min_{t \in \R^{p_{i_k}}}\Big\{\iprods{\nabla_i^k, t} + g_{i_k}(t + \tilde{z}^k _{i_k}) 
	                                           + \tfrac{\tau_k B^{k} _{i_{k}}}{2\tau_0} \Vert t \Vert_{(i_k)}^2 \Big\}.
	     \end{equation*}
           \State\label{step:z_tilde}
               Update $\tilde{z}^{k+1}_{i_k} := \tilde{z}^k _{i_k} + t^{k+1}_{i_{k}}$.
           \State\label{step:u}
               Update ${u}^{k+1}_{i_k} := {u}^k _{i_k} - \frac{1-\tau_k/\tau_0}{c_k} t^{k+1}_{i_{k}}$.

	   \State\label{step:update_tau_eff} Compute $\tau_{k+1} \in (0, 1)$ as the unique positive root of $\tau^3 + \tau^2 + \tau_k^2 \tau - \tau_k ^2=0$.
	   \State\label{step:update_beta_eff} Update $\beta_{k+2} := \frac{\beta_{k+1}}{1+\tau_{k+1}}$ and $B^{k+1} _i := \hat{L}_i + \frac{\Vert A_i \Vert ^2}{\beta_{k+2}}$ for $i\in [n]$.
     \EndFor
\end{algorithmic}
\end{algorithm}
\vspace{-3ex}
We present the following result which shows the equivalence between Algorithm~\ref{alg:A1} and Algorithm~\ref{alg:A_eff}, the proof of which can be found in the supplementary document.
\begin{proposition}\label{prop:eff_equiv}
Let $c_k = \prod_{l=0}^k (1 - \tau_l)$, $\hat{z}^k = c_k u^k + \tilde{z}^k$ and $\bar{z}^k = c_{k-1} u^k + \tilde{z}^k$. Then, $\tilde{x}^k=\tilde{z}^k$, $\hat{x}^k = \hat{z}^k$ and $\bar{x}^k=\bar{z}^k$, for all $k \geq 0$, where $\tilde{x}^k$, $\hat{x}^k$, and $\bar{x}^k$ are defined in Algorithm~\ref{alg:A1}.
\end{proposition}
\vspace{-1ex}
According to Algorithm ~\ref{alg:A_eff}, we never need to form or update full-dimensional vectors. Only times that we need $\hat{x}^k$ are when computing the gradient and the dual variable $y^{\ast}_{\beta_{k+1}}$. We present two special cases which are common in machine learning, in which we can compute these steps efficiently. 

\begin{remark}\label{rem:iter_complexity}
Under the following assumptions, we can characterize the per-iteration complexity explicitly. Let $A, M \in \mathbb{R}^{m\times p}$, and
\vspace{-1ex}
\begin{itemize}
\item[(a)] $f$ has the form $f(x) = \sum_{j=1}^m \phi_j(e_j^\top M x)$, where $e_j$ is the $j^{\text{th}}$ standard unit vector.
\vspace{-0.5ex}
\item[(b)] h is separable as in $h(Ax) = \delta_{\{c\}}(Ax)$ or $h(Ax) = \Vert Ax \Vert _1$.
\vspace{-1ex}
\end{itemize}

Assuming that we store and maintain the residuals $r_{u, f}^k = Mu^k$, $r_{\tilde{z}, f}^k = M\tilde{z}^k$, $r_{u, h}^k = Au^k$, $r_{\tilde{z}, h}^k = A\tilde{z}^k$, then we have the per-iteration cost as $\mathcal{O}(\max \{ | \{j \mid A_{ji} \neq 0\} |, |\{j \mid M_{ji} \neq 0\}|\})$ arithmetic operations.
If $h$ is partially separable as in~\cite{richtarik2016parallel}, then the complexity of each iteration will remain moderate.


\end{remark}

\vspace{-2ex}
\subsection{Case 1: Convergence analysis of SMART-CD for Lipschitz continuous $h$}\label{subsec:case1}
\vspace{-2ex}
We provide the following main theorem, which characterizes the convergence rate of Algorithm~\ref{alg:A1}.
\begin{theorem}\label{th:convergence1}
Let $x^{\star}$ be an optimal solution of \eqref{eq:cvx_prob} and let $\beta_1 > 0$ be given. 
In addition, let $\tau_0 := \min\set{ q_i \mid i\in[n]} \in (0, 1]$ and $\beta_0 := (1+\tau_0)\beta_1$ be given parameters.
For all $k\geq 1$, the sequence $\set{\bar{x}^k}$ generated by Algorithm \ref{alg:A1} satisfies:
\begin{equation}\label{eq:main_result1}
\expect{}{F(\bar{x}^{k}) - F^{\star}} \leq \frac{C^{\ast}(x^0)}{\tau_0(k-1) + 1} + \frac{\beta_1(1+\tau_0)D_{h^{\ast}}^2}{2(\tau_0k+1)},
\end{equation}
where $C^{\ast}(x^0) := (1-\tau_0)(F_{\beta_0}(x^0) - F^{\star}) + \sum_{i=1}^n{\!\!}\frac{\tau_0B^{0}_{i}}{2q_{i}} \Vert x^{\star}_{i} - x^{0}_{i} \Vert_{(i)}^2$ and $D_{h^{\ast}}$ is as defined by \eqref{eq:bound_est}.
\end{theorem}

%
In the special case when we use uniform distribution, $\tau_0 = q_i = 1/n$, the convergence rate reduces to 
\vspace{-.5ex}
\begin{align*}
\expect{}{F(\bar{x}^k) - F^\star} \leq \frac{nC^{\ast}(x^0)}{k+n-1} + \frac{(n+1)\beta_0 D_{h^\ast}^2}{2k+2n},
\end{align*}
where $C^{\ast}(x^0) := (1-\frac{1}{n})(F_{\beta_0}(x^0) - F^{\star}) + \sum_{i=1}^n{\!\!}\frac{B^{0}_{i}}{2} \Vert x^{\star}_{i} - x^{0}_{i} \Vert_{(i)}^2$.
This estimate shows that the convergence rate of Algorithm~\ref{alg:A1} is $$\mathcal{O}\left(\frac{n}{k}\right),$$ which is the best known so far to the best of our knowledge.


\vspace{-1.5ex}
\subsection{Case 2: Convergence analysis of SMART-CD for non-smooth constrained optimization}\label{sec:alg_constrained}
\vspace{-1.5ex}
In this section, we instantiate Algorithm~\ref{alg:A1} to solve constrained convex optimization problem with possibly non-smooth terms in the objective.
Clearly, if we choose $h(\cdot) = \delta_{\set{c}}(\cdot)$ in \eqref{eq:cvx_prob} as the indicator function of the set $\set{c}$ for a given vector $c\in\R^m$, then we obtain a constrained  problem:
\begin{equation}\label{eq:constr_cvx}
F^{\star} := \min_{x\in\R^p} \set{ F(x) = f(x) + g(x) \mid  Ax = c},
\end{equation}
where $f$ and $g$ are defined as in \eqref{eq:cvx_prob}, $A\in\R^{m\times p}$, and $c\in\R^m$.

We can specify Algorithm~\ref{alg:A1} to solve this constrained problem by modifying the following two steps:
\begin{itemize}
\vspace{-1.5ex}
\item[(a)] The update of $y^{\ast}_{\beta_{k+1}}(A \hat{x}^k)$ at Step~\ref{step:y_hat} is changed to 
\begin{equation}\label{eq:update_yhat2}
y^{\ast}_{\beta_{k+1}}(A\hat{x}^k) := \dot{y} + \tfrac{1}{\beta_{k+1}}(A\hat{x}^k - c),
\end{equation}
which  requires one matrix-vector multiplication in $A\hat{x}^k$. 
\item[(b)] The update of $\tau_k$ at Step~\ref{step:update_tau} and $\beta_{k+1}$ at Step~\ref{step:update_beta} are changed to
\begin{equation}\label{eq:update_tau_beta}
\tau_{k+1} := \tfrac{\tau_k}{1 + \tau_k}~~~\text{and}~~\beta_{k+2} := (1-\tau_{k+1})\beta_{k+1}.
\end{equation}
\end{itemize}
Now, we analyze the convergence of this algorithm by providing the following theorem.
\begin{theorem}\label{th:convergence2}
Let $\set{\bar{x}^k}$ be the sequence generated by Algorithm~\ref{alg:A1} for solving \eqref{eq:constr_cvx} using the updates \eqref{eq:update_yhat2} and \eqref{eq:update_tau_beta} and let $y^\star$ be an arbitrary optimal solution of the dual problem of~\eqref{eq:constr_cvx}. In addition, let $\tau_0 := \min\set{ q_i \mid i\in[n]} \in (0, 1]$ and $\beta_0 := (1+\tau_0)\beta_1$ be given parameters. Then, we have the following estimates:
\begin{equation}\label{eq:main_result2}
\left\{\begin{array}{lll}
&\expect{}{F(\bar{x}^k) - F^{\star}} &\leq  \frac{C^{\ast}(x^0)}{\tau_0(k-1) + 1} + \frac{\beta_1\Vert y^\star - \dot y\Vert^2}{2(\tau_0(k-1)+1)} + \Vert y^{\star}\Vert\expect{}{\Vert A\bar{x}^k - b\Vert}, \vspace{1.5ex}\\
&\expect{}{\Vert A\bar{x}^k - b\Vert} &\leq \frac{\beta_1}{\tau_0(k-1) + 1}\left[\Vert y^{\star}-\dot{y}\Vert + \left(\Vert y^{\star}-\dot{y}\Vert^2 + 2\beta_1^{-1}C^{\ast}(x^0)\right)^{1/2}\right],
\end{array}\right.
\end{equation}
where $C^{\ast}(x^0) := (1-\tau_0)(F_{\beta_0}(x^0) - F^{\star}) + \sum_{i=1}^n{\!\!}\frac{\tau_0B^{0}_{i}}{2q_{i}} \Vert x^{\star}_{i} - x^{0}_{i} \Vert_{(i)}^2$.
We note that the following lower bound always holds $-\Vert y^{\star}\Vert\expect{}{\Vert A\bar{x}^k - b\Vert} \leq \expect{}{F(\bar{x}^k) - F^{\star}}$.
\end{theorem}

\vspace{-2ex}
\subsection{Other special cases}\label{subsec:spec_cases}
\vspace{-2ex}
We consider the following special cases of Algorithm~\ref{alg:A1}: 

\vspace{-2ex}
\paragraph{The case $h = 0$:} 
In this case, we obtain an algorithm similar to the one studied in \cite{qu2016coordinate} except that we have non-uniform sampling instead of importance sampling. If the distribution is uniform, then we obtain the method in \cite{fercoq2015accelerated}.

\vspace{-2ex}
\paragraph{The case $g = 0$:} 
In this case, we have $F(x) = f(x) + h(Ax)$, which can handle the linearly constrained problems with smooth objective function.
In this case, we can choose $\tau_0 = 1$, and the coordinate proximal gradient step, Step~\ref{step:x_tilde} in Algorithm~\ref{alg:A1}, is simplified as
\begin{equation}\label{eq:coordinate_grad_step}
   \tilde{x}^{k+1}_{i_{k}} := \tilde{x}^k_{i_k} - \tfrac{q_{i_k}}{\tau_kB_{i_k}^k}H_{i_k}^{-1}\left( \nabla_{i_k}f(\hat{x}^k) + A_{i_k}^\top y^{\ast}_{\beta_{k+1}}(\hat{u}^k)\right).
\end{equation}
In addition, we replace Step~\ref{step:x_bar} with 
\vspace{-1ex}
\begin{equation}\label{eq:st8_ch}
\bar{x}^{k+1}_i = \hat{x}^k_i + \frac{\tau_k}{q_i}(\tilde{x}^{k+1}_i - \tilde{x}^k_i), ~~ \forall i \in [n].
\end{equation}
\vspace{-2ex}
We then obtain the following results:
\begin{corollary}\label{co:zero_g_case}
Assume that Assumption~\ref{as:fun_assumption} holds.
Let $\tau_0 = 1$, $\beta_1 > 0$ and Step~\ref{step:x_tilde} and~\ref{step:x_bar} of Algorithm~\ref{alg:A1} be updated by \eqref{eq:coordinate_grad_step} and \eqref{eq:st8_ch}, respectively. 
If, in addition, $h$ is Lipschitz continuous, then we have
\begin{equation}\label{eq:main_result1c}
\expect{}{F(\bar{x}^{k}) - F^{\star}} \leq \frac{1}{k}\sum_{i=1}^n{\!\!}\frac{B^{0}_{i}}{2q^2_{i}} \Vert x^{\star}_{i} - x^{0}_{i} \Vert_{(i)}^2 + \frac{\beta_1D_{h^{\ast}}^2}{k+1},
\end{equation}
where $D_{h^{\ast}}$ is defined by \eqref{eq:bound_est}.

If, instead of Lipschitz continuous $h$, we have $h(\cdot) = \delta_{\set{c}}(\cdot)$ to solve the constrained problem~\eqref{eq:constr_cvx} with $g = 0$, then we have
\begin{equation}\label{eq:main_result2c}
\left\{\begin{array}{lll}
&\expect{}{F(\bar{x}^k) - F^{\star}} &\leq  \frac{C^{\ast}(x^0)}{k} + \frac{\beta_1 \Vert y^\star - \dot y\Vert^2}{2k} + \Vert y^{\star}\Vert\expect{}{\Vert A\bar{x}^k - b\Vert}, \vspace{1.5ex}\\
&\expect{}{\Vert A\bar{x}^k - b\Vert} &\leq \frac{\beta_1}{k}\left[\Vert y^{\star}-\dot{y}\Vert + \left(\Vert y^{\star}-\dot{y}\Vert^2 + 2\beta_1^{-1}C^{\ast}(x^0)\right)^{1/2}\right],
\end{array}\right.
\end{equation}
where $C^{\ast}(x^0) := \sum\limits_{i=1}^n{\!\!}\frac{B^{0}_{i}}{2q^2_{i}} \Vert x^{\star}_{i} - x^{0}_{i} \Vert_{(i)}^2$.
\end{corollary}
%

\vspace{-2ex}
\subsection{Restarting SMART-CD}
\vspace{-2ex}
It is known that restarting an accelerated method significantly enhances its practical performance when the underlying problem admits a (restricted) strong convexity condition. 
As a result, we describe below how to restart (i.e., the momentum term) in Efficient SMART-CD.
If the restart is injected in the $k$-th iteration, then we restart the algorithm with the following steps:
\vspace{-.5ex}
\begin{equation*}
\left\{\begin{array}{lll}
u^{k+1} &\leftarrow 0, \\
r^{k+1}_{u, f} &\leftarrow 0, \\ 
r^{k+1}_{u, h} &\leftarrow 0, \\
\dot{y} &\leftarrow y_{\beta_{k+1}}^{\ast}(c_k r^k_{u, h} + r^k_{\tilde{z}, h}), \\
\beta_{k+1} &\leftarrow \beta_1, \\
\tau_{k+1} &\leftarrow \tau_0, \\
c_k &\leftarrow 1.
\end{array}\right.
\end{equation*}
\vspace{-2.5ex}

The first three steps of the restart procedure is for restarting the primal variable which is classical \cite{o2015adaptive}. 
Restarting $\dot{y}$ is also suggested in \cite{tran2015smooth}. The cost of this procedure is essentially equal to the cost of one iteration as described in Remark~\ref{rem:iter_complexity}, therefore even restarting once every epoch will not cause a significant difference in terms of per-iteration cost.

\vspace{-2ex}
\section{Numerical evidence}\label{sec:numerical_experiments}
\vspace{-2ex}
We illustrate the performance of Efficient SMART-CD in brain imaging and support vector machines applications. We also include one representative example of a degenerate linear program to illustrate why the convergence rate guarantees of our algorithm matter. We compare SMART-CD with Vu-Condat-CD \cite{fercoq2015coordinate}, which is a coordinate descent variant of Vu-Condat's algorithm \cite{vu2013splitting}, FISTA \cite{beck2009fast}, ASGARD \cite{tran2015smooth}, Chambolle-Pock's primal-dual algorithm \cite{chambolle2011first}, L-BFGS \cite{byrd1995limited} and SDCA \cite{shalev2013stochastic}.

\vspace{-2ex}
\subsection{A degenerate linear program: Why do convergence rate guarantees matter? }\label{sec:lp_exp}
\vspace{-2ex}
We consider the following degenerate linear program studied in \cite{tran2015smooth}:
\begin{equation}\label{eq:lp_statement}
\left\{\begin{array}{lll}
&\displaystyle\min_{x \in \R^p} & 2 x_p \\
& \text{s.t.}  &\sum_{k=1}^{p-1} x_k = 1, \\
& & x_p - \sum_{k=1}^{p-1} x_k = 0, \qquad (2 \leq j \leq d), \\
& & x_p \geq 0.
\end{array}\right.
\end{equation}
Here, the constraint $x_p - \sum_{k=1}^{p-1} x_k = 0$ is repeated $d$ times.
This problem satisfies the linear constraint qualification condition, which guarantees the primal-dual optimality. 
If we define
\begin{equation*}
f(x) = 2x_p, \quad g(x) = \delta_{\{x_p \geq 0\}}(x_p), \quad h(Ax) = \delta_{\{c\}}(Ax),
\end{equation*}
where
\vspace{-1ex}
\begin{equation*}
Ax = \left[\sum_{k=1}^{p-1}x_k,~~x_p - \sum_{k=1}^{p-1}x_k,\dots,~~ x_p - \sum_{k=1}^{p-1}x_k \right]^\top, \quad c=[1, 0, \dots, 0]^\top,
\end{equation*}
we can fit this problem and its dual form into our template \eqref{eq:cvx_prob}.

For this experiment, we select the dimensions $p=10$ and $d=200$.
We implement our algorithm and compare it with Vu-Condat-CD. We also combine our method with the restarting strategy proposed above.
We use the same mapping to fit the problem into the template of Vu-Condat-CD. 

\begin{figure}[ht!]
\centering
\begin{tabular}{ccc}
\hspace{-2.5mm}\includegraphics[width=0.45\columnwidth]{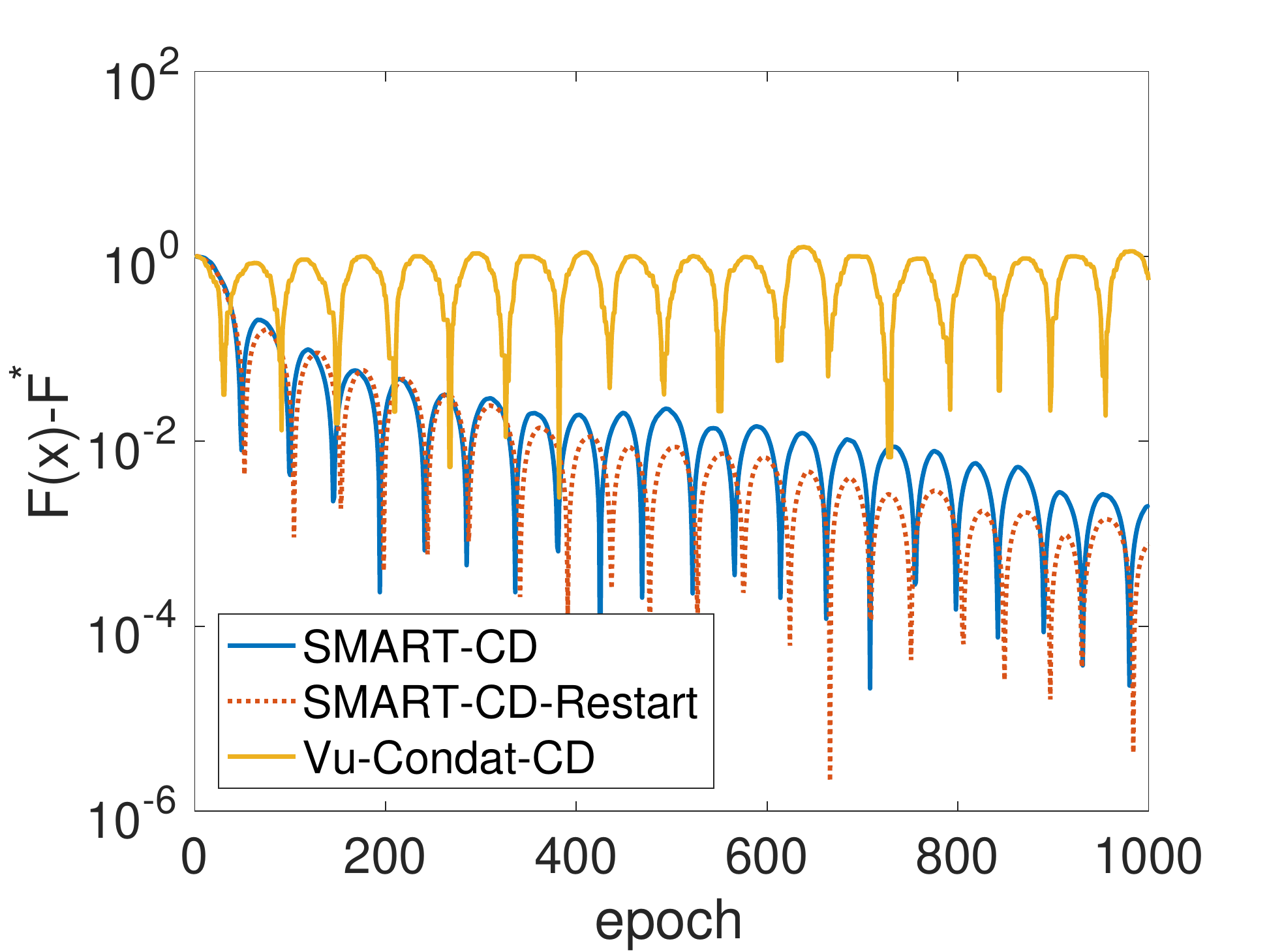} &
\hspace{-2.5mm}\includegraphics[width=0.45\columnwidth]{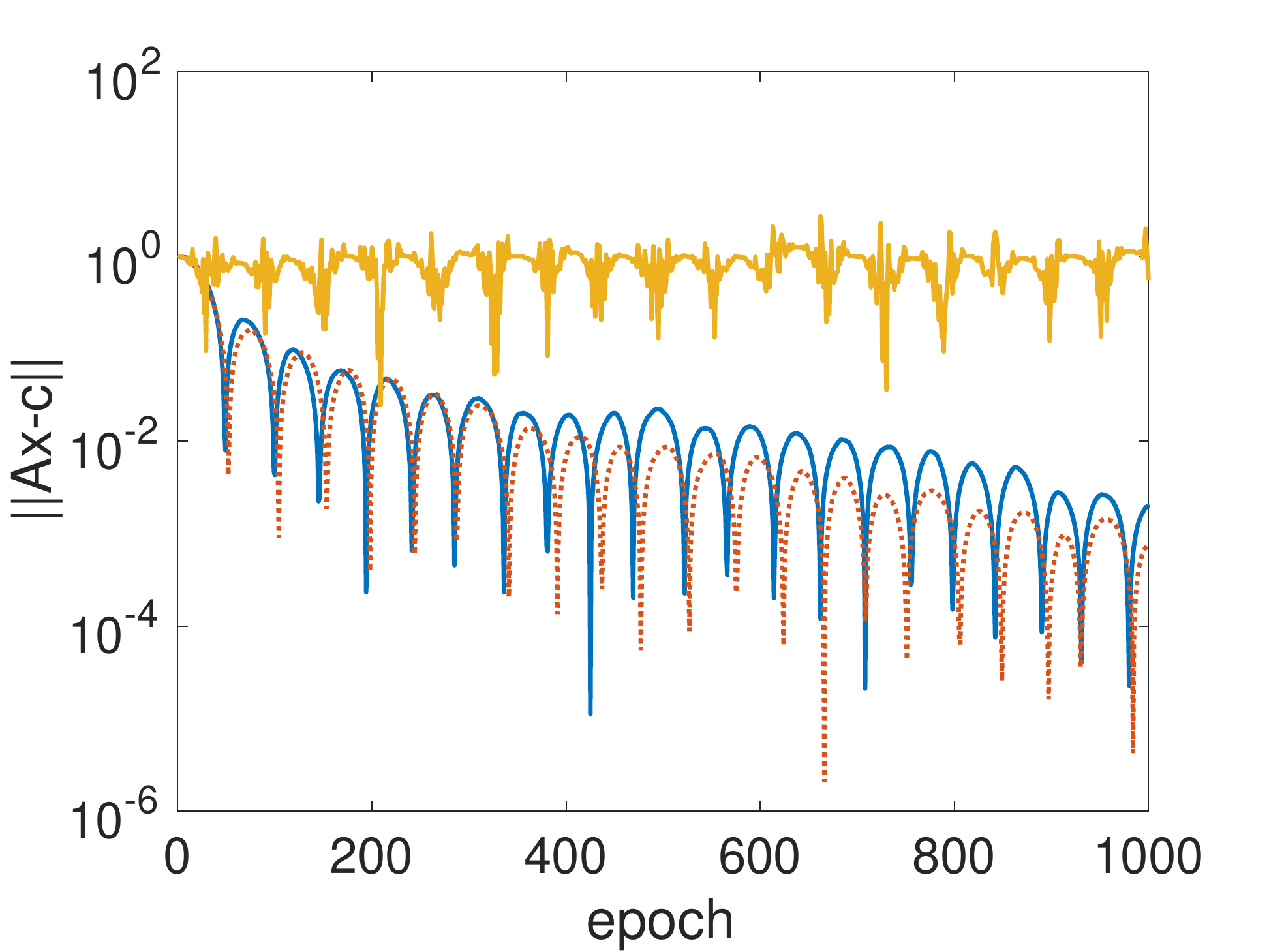} 
\vspace{-1ex}
\end{tabular}
\vspace{-1ex}
\caption{The convergence behavior of $3$ algorithms on a degenerate linear program.}\label{fig:lp_exam}
\vspace{-2ex}
\end{figure}

Figure~\ref{fig:lp_exam} illustrates the convergence behavior of Vu-Condat-CD and SMART-CD. We compare primal suboptimality and feasibility in the plots. The explicit solution of the problem is used to generate the plot with primal suboptimality. We observe that degeneracy of the problem prevents Vu-Condat-CD from making any progress towards the solution, where SMART-CD preserves $\mathcal{O}(1/k)$ rate as predicted by theory. 
We emphasize that the authors in \cite{fercoq2015coordinate} proved almost sure convergence for Vu-Condat-CD but they did not provide a convergence rate guarantee for this method.
Since the problem is certainly non-strongly convex, restarting does not significantly improve performance of SMART-CD.

\vspace{-2ex}
\subsection{Total Variation and $\ell_1$-regularized least squares regression with functional MRI data}\label{sec:tv_exp}
\vspace{-2ex}
In this experiment, we consider a computational neuroscience application where prediction is done based on a sequence of functional MRI images. Since the images are high dimensional and the number of samples that can be taken is limited, TV-$\ell_1$ regularization is used to get stable and predictive estimation results \cite{dohmatob2014benchmarking}. The convex optimization problem we solve is of the form:
\begin{equation}\label{eq:tv_problem}
\min_{x \in \R^p} \tfrac{1}{2} \Vert Mx - b \Vert ^2 + \lambda r\Vert x \Vert _1 + \lambda (1-r) \Vert x \Vert_{\mathrm{TV}}.
\end{equation}
This problem fits to our template with
\begin{equation*}
f(x) = \tfrac{1}{2} \Vert Mx - b \Vert ^2, \qquad 
g(x) = \lambda r \Vert x \Vert _1, \qquad 
h(u) = \lambda (1-r) \Vert u \Vert_{1},
\end{equation*}
where $D$ is the 3D finite difference operator to define a total variation norm $\Vert\cdot\Vert_{\mathrm{TV}}$ and $u = Dx$.

We use an fMRI dataset where the primal variable $x$ is 3D image of the brain that contains $33177$ voxels. 
Feature matrix $M$ has $768$ rows, each representing the brain activity for the corresponding example \cite{dohmatob2014benchmarking}. 
We compare our algorithm with Vu-Condat's algorithm, FISTA, ASGARD, Chambolle-Pock's primal-dual algorithm, L-BFGS and Vu-Condat-CD. 

\begin{figure}[ht!]
\vspace{-2ex}
\centering
\begin{tabular}{ccc}
\hspace{-2.5mm}\includegraphics[width=0.34\columnwidth]{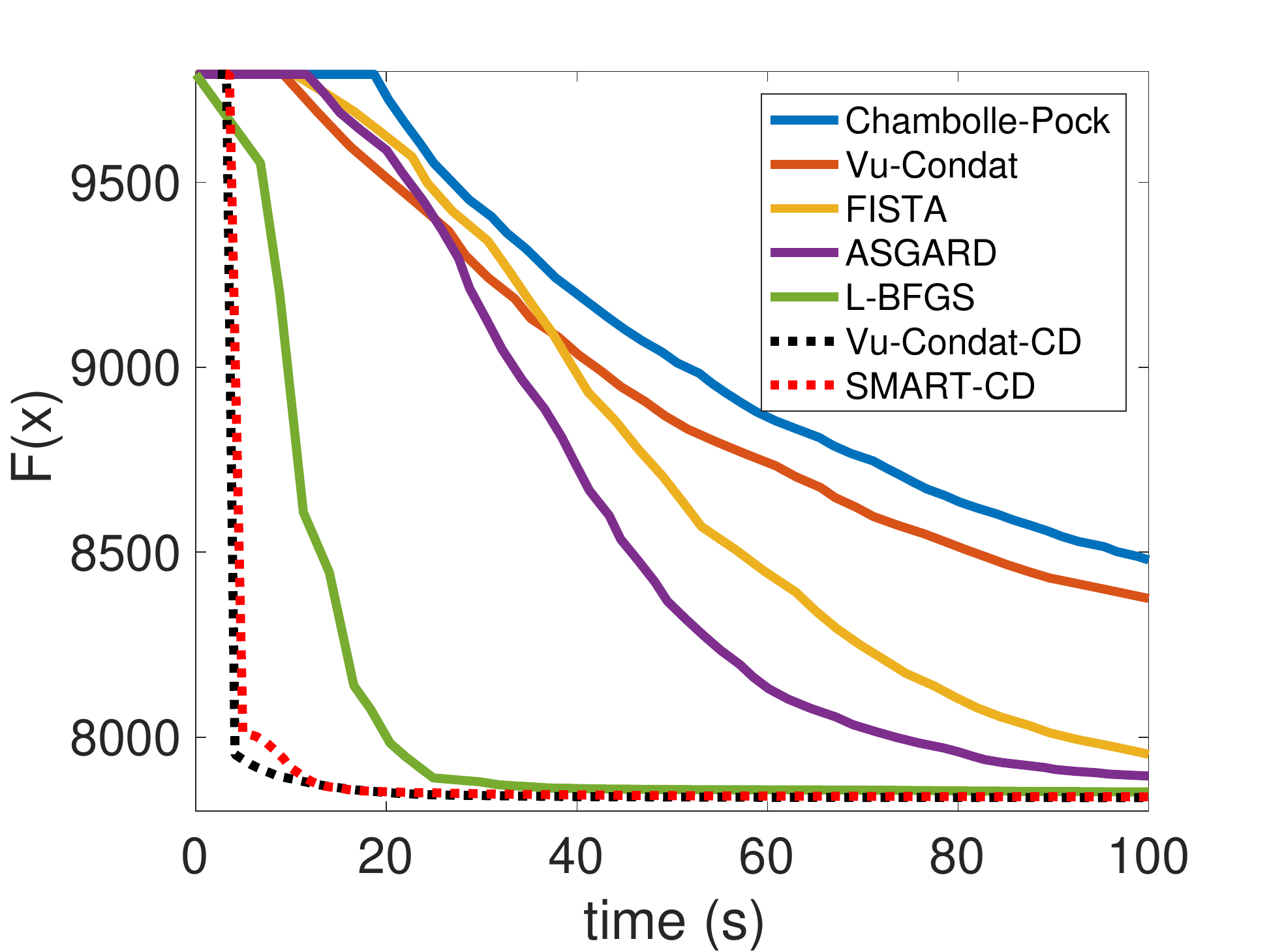}
\hspace{-2.5mm}\includegraphics[width=0.34\columnwidth]{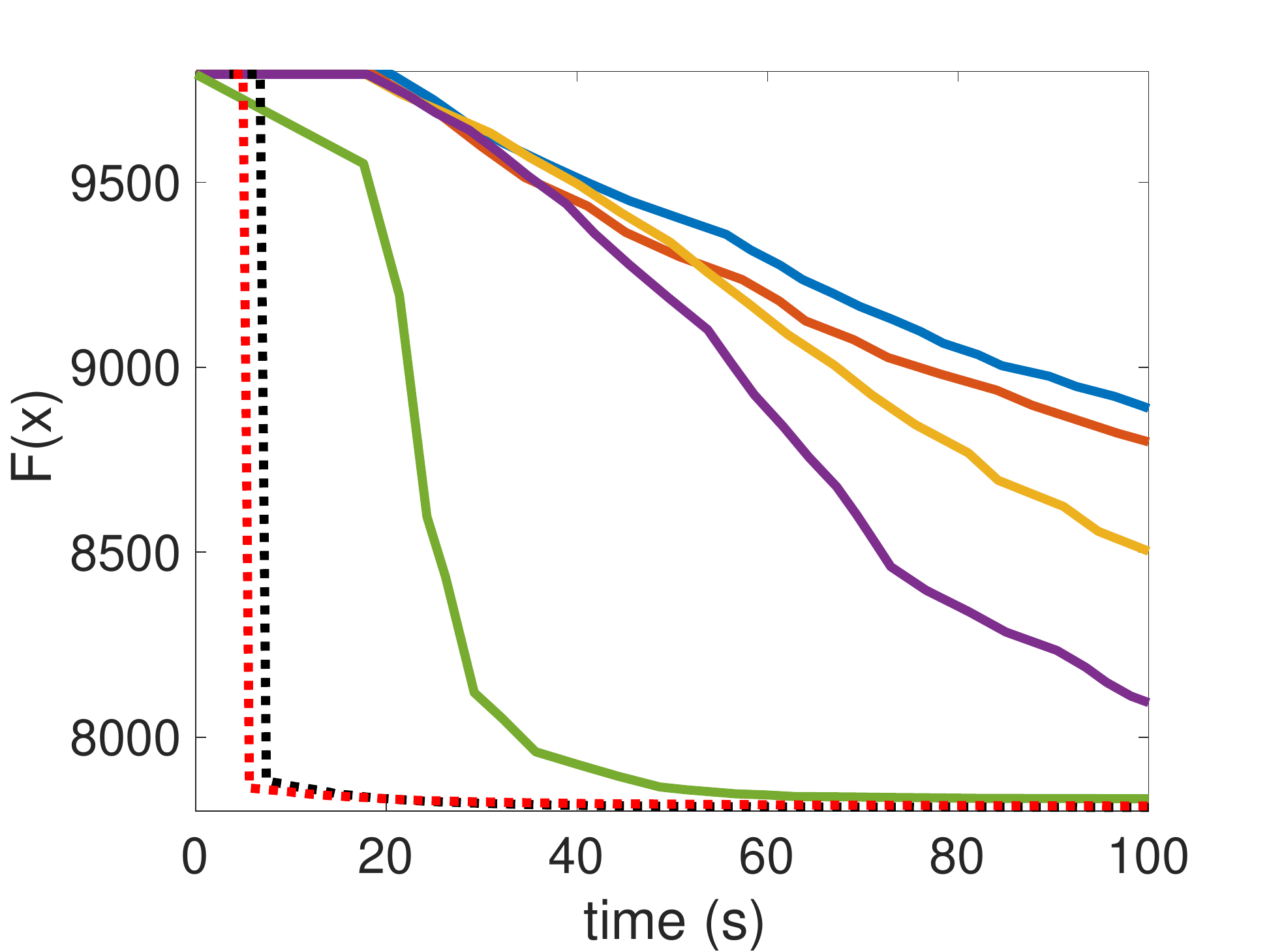} &
\hspace{-2.5mm}\includegraphics[width=0.34\columnwidth]{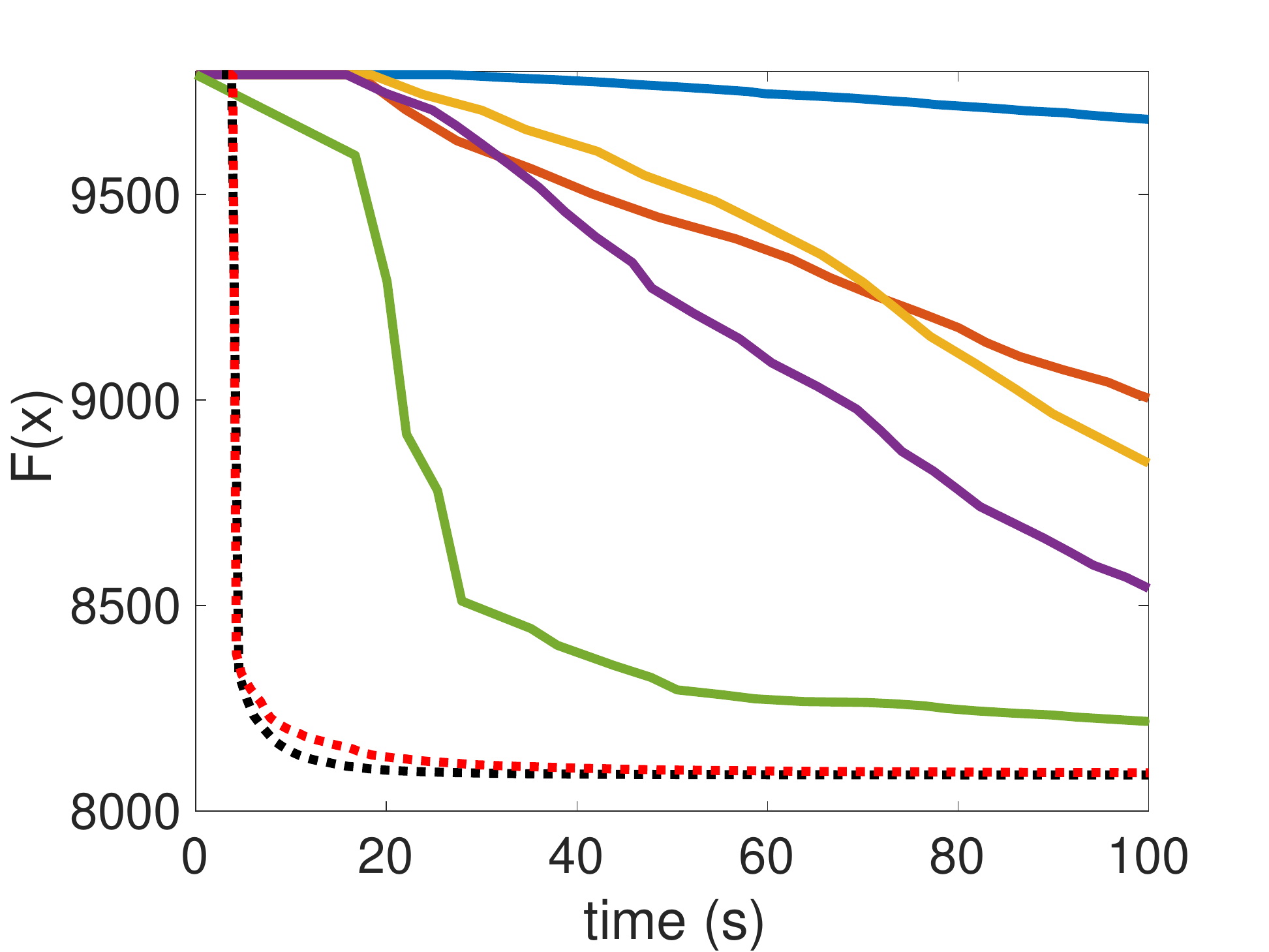} 
\end{tabular}
\vspace{-2ex}
\caption{The convergence of $7$ algorithms for problem~\eqref{eq:tv_problem}. 
The regularization parameters for the first plot are $\lambda = 0.001, r=0.5$, for the second plot are $\lambda = 0.001, r=0.9$, for the third plot are $\lambda = 0.01, r=0.5$ .}\label{fig:tv_exam}
\vspace{-2ex}
\end{figure}
Figure~\ref{fig:tv_exam} illustrates the convergence behaviour of the algorithms for different values of the regularization parameters. 
Per-iteration cost of SMART-CD and Vu-Condat-CD is similar, therefore the behavior of these two algorithms are quite similar in this experiment. 
Since Vu-Condat's, Chambolle-Pock's, FISTA and ASGARD methods work with full dimensional variables, they have slow convergence in time. L-BFGS has a close performance to coordinate descent methods.

The simulation in Figure \ref{fig:tv_exam} is performed using benchmarking tool of \cite{dohmatob2014benchmarking}. The algorithms are tuned for the best parameters in practice.

\vspace{-2ex}
\subsection{Linear support vector machines problem with bias}\label{sec:dual_svm_exp}
\vspace{-2ex}
In this section, we consider an application of our algorithm to support vector machines (SVM) problem for binary classification. 
Given a training set with $m$ examples $\left\{ a_1, a_2, \dots, a_m \right\}$ such that $a_i \in \R^p$ and class labels $\left\{ b_1, b_2, \dots\, b_m \right\}$ such that $b_i \in \{ -1, +1 \}$, we define the soft margin primal support vector machines problem with bias as
\begin{equation}\label{eq:primal_svm_problem}
\min_{w \in \R^{p}} \sum_{i=1}^{m} C_i \max \Big(0, 1 - b_i (\iprods{a_i, w} + w_0 ) \Big) + \tfrac{\lambda}{2} \Vert w \Vert ^2.
\end{equation}
As it is a common practice, we solve its dual formulation, which is a constrained problem:
\begin{equation}\label{eq:dual_svm_problem}
\left\{\begin{array}{lll}
&\min\limits_{x\in\R^{m}} &\set{ \frac{1}{2\lambda} \Vert M D(b) x \Vert ^2 - \sum_{i=1}^m x_i } \vspace{1.5ex}\\
& \mathrm{s.t.} & 0 \leq x_i \leq C_i,  ~~i=1,\cdots, m,~~~b^{\top}x=0,
\end{array}\right.
\end{equation}
where $D(b)$ represents a diagonal matrix that has the class labels $b_i$ in its diagonal and $M \in \R^{p\times m}$ is formed by the example vectors. 
If we define
\vspace{-0.5ex}
\begin{equation*}
f(x) = \frac{1}{2\lambda} \Vert M D(b) x \Vert ^2 - \sum_{i=1}^{m} x_i, \quad 
g_i(x_i) = \delta_{\{ 0 \leq x_i \leq C_i \}}, \quad 
c=0, \quad 
A = b^{\top},
\end{equation*}
then, we can fit this problem into our template in \eqref{eq:constr_cvx}.

\begin{figure}[ht!]
\vspace{-2ex}
\centering
\begin{tabular}{ccc}
\hspace{-2.5mm}\includegraphics[width=0.34\columnwidth]{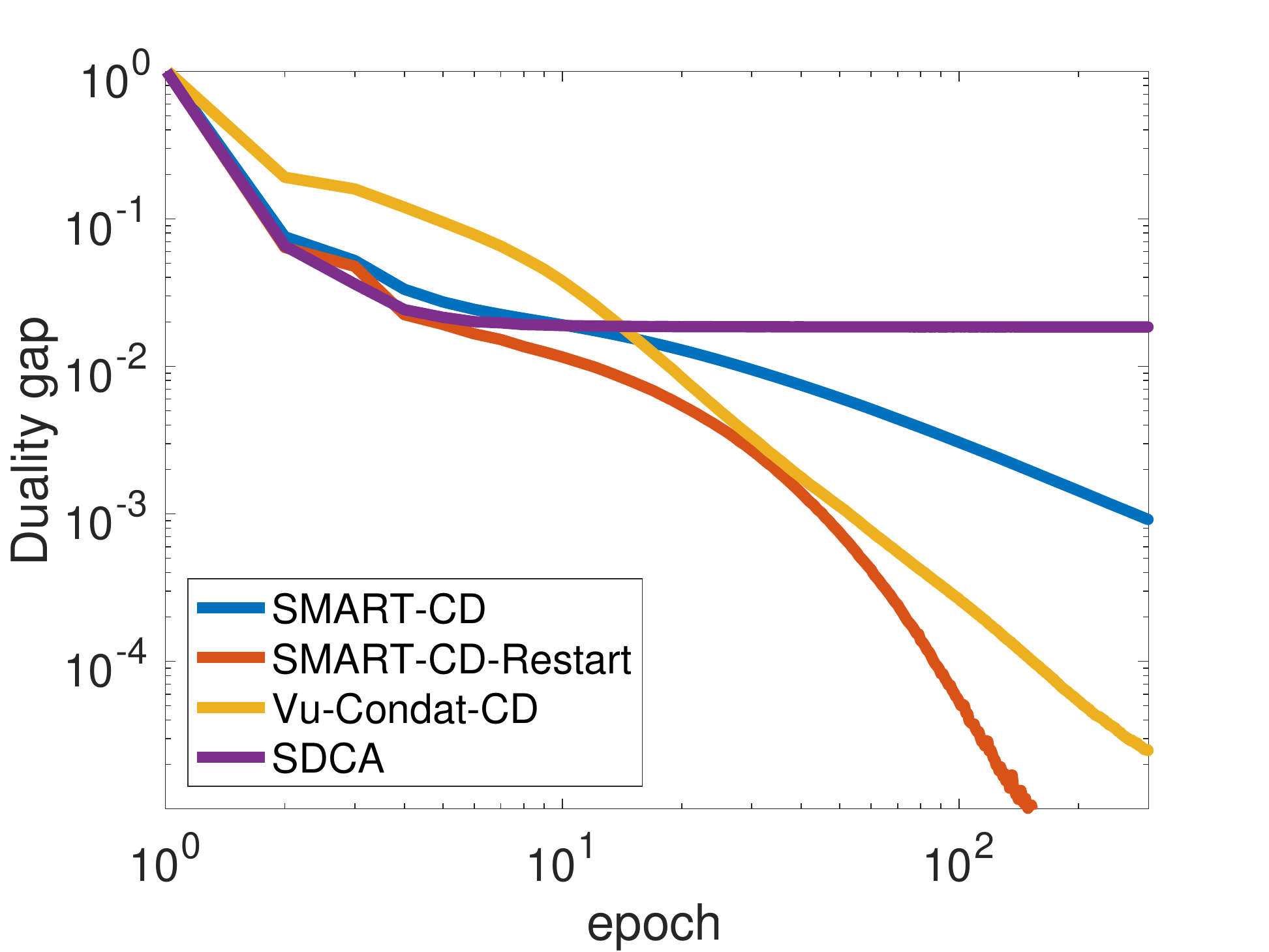}
\hspace{-2.5mm}\includegraphics[width=0.34\columnwidth]{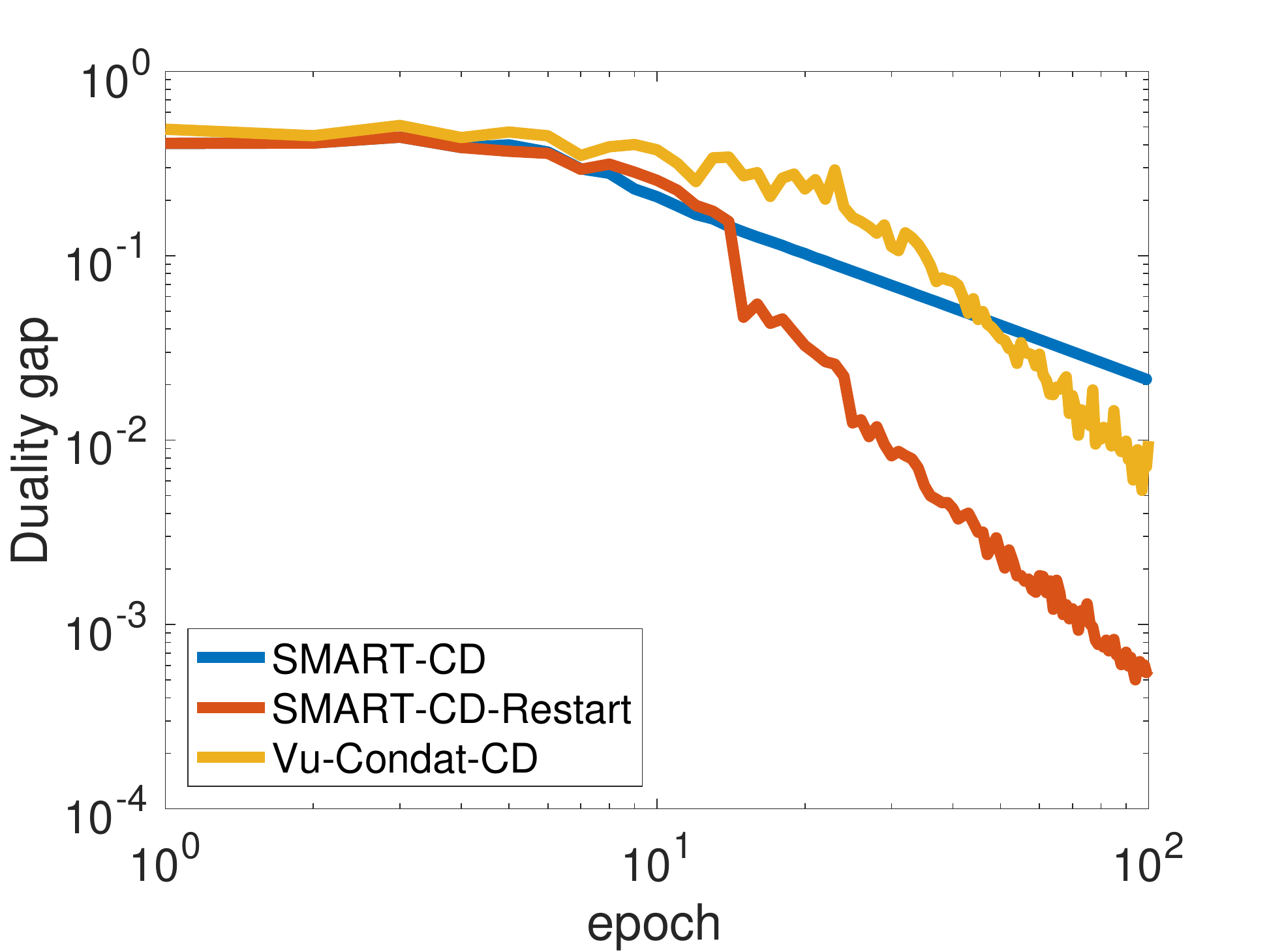} &
\hspace{-2.5mm}\includegraphics[width=0.34\columnwidth]{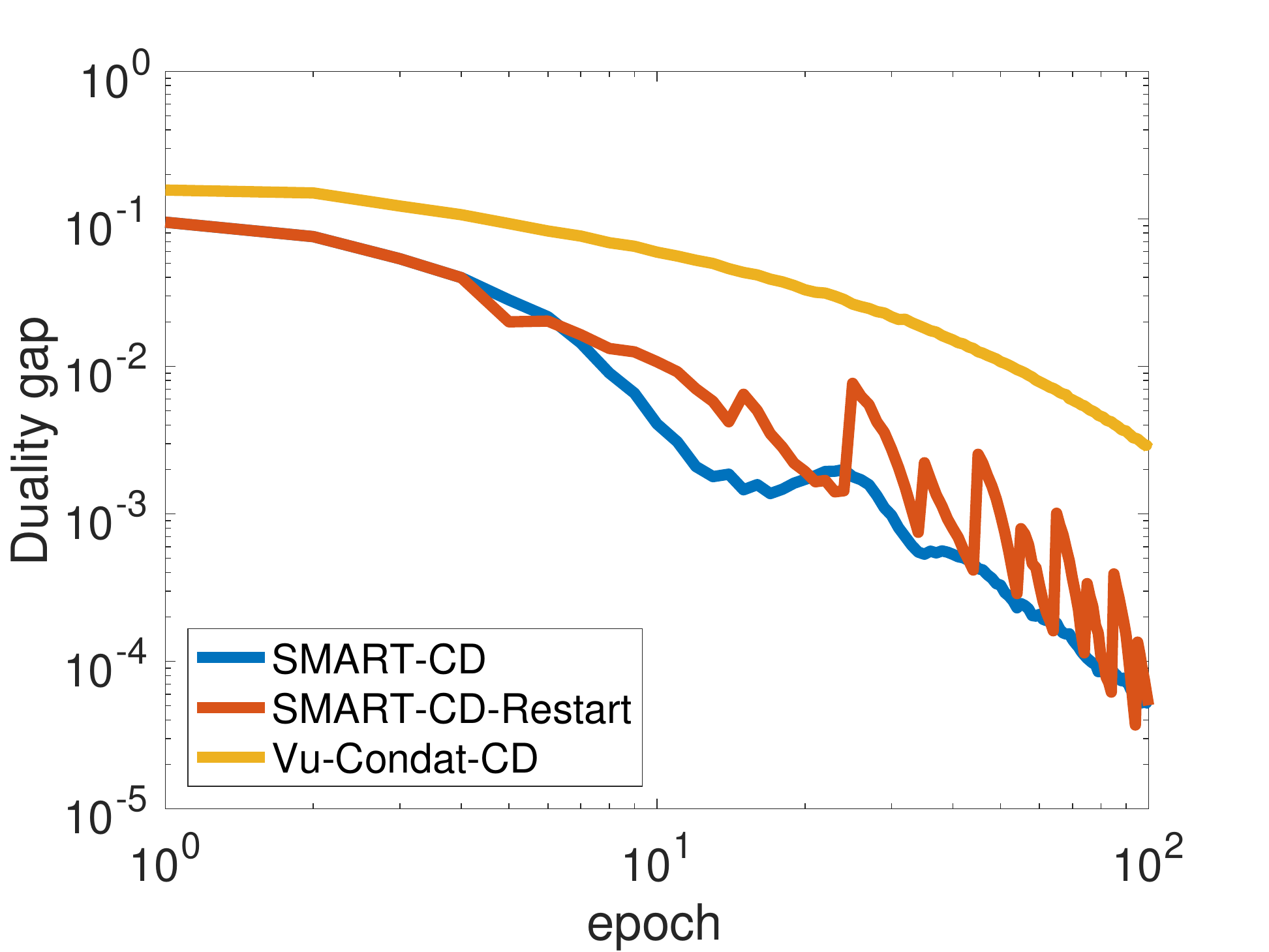} 
\end{tabular}
\vspace{-2ex}
\caption{The convergence of $4$ algorithms on the dual SVM \eqref{eq:dual_svm_problem} with bias. We only used SDCA in the first dataset since it stagnates at a very low accuracy.}\label{fig:rcv1_plot}
\vspace{-1ex}
\end{figure}

\vspace{-0.5ex}
We apply the specific version of SMART-CD for constrained setting from Section~\ref{sec:alg_constrained} and compare with Vu-Condat-CD and SDCA. 
Even though SDCA is a state-of-the-art method for SVMs, we are not able to handle the bias term using SDCA. 
Hence, it only applies to~\eqref{eq:dual_svm_problem} when $b^{\top}x=0$ constraint is removed. 
This causes SDCA not to converge to the optimal solution when there is bias term in the problem~\eqref{eq:primal_svm_problem}. The following table summarizes the properties of the classification datasets we used.
\vspace{-1ex}
\begin{center}
    \begin{tabular}{| l | l | l | l |}
    \hline
    Data Set & Training Size & Number of Features & Convergence Plot \\ \hline
    rcv1.binary \cite{chang2011libsvm, lewis2004rcv1} & 20,242 & 47,236 & Figure~\ref{fig:rcv1_plot}, plot 1 \\ \hline
    a8a \cite{chang2011libsvm, Lichman:2013} & 22,696 & 123 & Figure~\ref{fig:rcv1_plot}, plot 2  \\ \hline
    gisette \cite{chang2011libsvm,guyon2005result} & 6,000 & 5,000 & Figure~\ref{fig:rcv1_plot}, plot 3 \\
    \hline
    \end{tabular}
\end{center}
\vspace{-1ex}
Figure~\ref{fig:rcv1_plot} illustrates the performance of the algorithms for solving the dual formulation of SVM in \eqref{eq:dual_svm_problem}. We compute the duality gap for each algorithm and present the results with epochs in the horizontal axis since per-iteration complexity of the algorithms is similar. 
As expected, SDCA gets stuck at a low accuracy since it ignores one of the constraints in the problem. We demonstrate this fact in the first experiment and then limit the comparison to SMART-CD and Vu-Condat-CD.
Equipped with restart strategy, SMART-CD shows the fastest convergence behavior due to the restricted strong convexity of~\eqref{eq:dual_svm_problem}.

\vspace{-2ex}
\section{Conclusions}
\vspace{-2ex}
Coordinate descent methods have been increasingly deployed to tackle huge scale machine learning problems in recent years. The most notable works include \cite{Nesterov2012,richtarik2014iteration,richtarik2016parallel,fercoq2015accelerated,shalev2013stochastic,necoara2016parallel}.
Our method relates to several works in the literature including \cite{fercoq2015accelerated,fercoq2013smooth,Nesterov2012,nesterov2017efficiency,qu2016coordinate,tran2015smooth}. 
The algorithms developed in \cite{fercoq2015accelerated,richtarik2014iteration,richtarik2016parallel} only considered a special case of \eqref{eq:cvx_prob} with $h = 0$, and cannot be trivially extended to apply to general setting \eqref{eq:cvx_prob}.
Here, our algorithm can be viewed as an adaptive variant of the method developed in \cite{fercoq2015accelerated} extended to the sum of three functions. 
The idea of homotopy strategies relate to \cite{tran2015smooth} for first-order primal-dual methods. 
This paper further extends such an idea to randomized coordinate descent methods for solving \eqref{eq:cvx_prob}.
We note that a naive application of the method developed in \cite{fercoq2015accelerated} to the smoothed problem with a carefully chosen fixed smoothness parameter would result in the complexity $\mathcal{O}(n^2/k)$, whereas using our homotopy strategy on the smoothness parameter, we reduced this complexity to $\mathcal{O}
(n/k)$.

With additional strong convexity assumption on problem template~\eqref{eq:cvx_prob}, it is possible to obtain $\mathcal{O}(1/k^2)$ rate by using deterministic first-order primal-dual algorithms \cite{chambolle2011first, tran2015smooth}.
It remains as future work to incorporate strong convexity to coordinate descent methods for solving nonsmooth optimization problems with a faster convergence rate.

\vspace{-2ex}
\section*{Acknowledgments}
\vspace{-2ex}
The work of O. Fercoq was supported by a public grant as part of the Investissement d'avenir project, reference ANR-11-LABX-0056-LMH, LabEx LMH. 
The work of Q. Tran-Dinh was partly supported by NSF grant, DMS-1619884, USA. The work of A. Alacaoglu and V. Cevher was supported by European Research Council (ERC) under the European Union's Horizon 2020 research and innovation programme (grant agreement n$^\text{o}$ 725594 - time-data).


\bibliographystyle{ieeetr}


\newpage
\appendix 
\begin{center}
\textbf{\large Supplementary document}\vspace{1ex}\\
\textbf{\Large Smooth Primal-Dual Coordinate Descent Algorithms for Nonsmooth Convex Optimization}
\end{center}

\vspace{-2ex}
\section{Key lemmas}
\vspace{-2ex}
The following properties are key to design the algorithm, whose proofs are very similar to the proof of \cite[Lemma 10]{tran2015smooth} by using a different norm, and we omit the proof here. The proof of the last property directly follows by using the explicit form of $h_\beta(u)$ in the special case when $h^\ast(y) = \langle c, y \rangle$.
\begin{lemma}\label{le:smoothing_properties}
For any $u, \hat{u}\in\R^m$, the function $h_{\beta}$ defined by \eqref{eq:smooth_h} satisfies the following properties:
\begin{itemize}
\item[$\mathrm{(a)}$] $h_{\beta}(\cdot)$ is convex and smooth. Its gradient $\nabla{h_{\beta}}(u) = y^{\ast}_{\beta}(u)$ is Lipschitz continuous with the Lipschitz constant $L_{h_{\beta}} = \frac{1}{\beta}$.
\item[$\mathrm{(b)}$] $h_{\beta}(u) + \iprods{\nabla{h_{\beta}}(u), \hat{u} - u} + \frac{\beta}{2}\Vert y^{\ast}_{\beta}(u) - y^{\ast}_{\beta}(\hat{u})\Vert^2 \leq h_{\beta}(\hat{u})$.
\item[$\mathrm{(c)}$] $h(\hat{u}) \geq h_{\beta}(u) +  \iprods{\nabla{h_{\beta}}(u), \hat{u} - u}  + \frac{\beta}{2}\Vert y^{\ast}_{\beta}(u) - \dot{y}\Vert^2$.
\item[$\mathrm{(d)}$] $h_{\beta}(u) \leq h_{\bar{\beta}}(u) + \left(\frac{\bar{\beta} - \beta}{2}\right)\Vert y_{\beta}^{\ast}(u) - \dot{y}\Vert^2$.
\item[$\mathrm{(e)}$] \label{it:lin} If $h^{\ast}(y) = \iprods{c,y}$, a linear function, then $h_{\beta}(u) = h_{\bar{\beta}}(u) + \frac{(\bar{\beta} - \beta)\beta}{2\bar{\beta}}\Vert y_{\beta}^{\ast}(u) - \dot{y}\Vert^2$.
\end{itemize}
\end{lemma}

\begin{lemma}\label{le:parameters}
The parameters $\{ \tau_k \}_{k\geq 0}$ and $\{\beta_k \}_{k\geq1}$ updated by Steps~\ref{step:update_tau} and~\ref{step:update_beta}, respectively, satisfy the following bounds:
\begin{align}\label{eq:bound_tau_beta}
\frac{1}{k+\tau_0^{-1}} \leq \tau_k \leq \frac{2}{k+\tau_0^{-1} + 1}, ~~~~~~~~~ \beta_k \leq \frac{\beta_1 (1+\tau_0)}{\tau_0 k + 1}.
\end{align}
\begin{proof}
We proceed by induction. By Step~\ref{step:update_tau}, we have $\tau_{k-1}^2 = \frac{\tau_k^3 + \tau_k^2}{1-\tau_k}$. For $k=0$, the bounds trivially hold since $\tau_0 \leq \frac{1}{n}$. By the inductive assumption, we have $\frac{1}{k-1+\tau_0^{-1}} \leq \tau_{k-1} \leq \frac{2}{k+\tau_0^{-1}}$. Assume toward condtradiction that $\tau_k < \frac{1}{k+\tau_0^{-1}}$. Then $\frac{1}{(k-1+\tau_0^{-1})^2} \leq \tau_{k-1}^2 = \frac{\tau_k^3 + \tau_k^2}{1-\tau_k} < \frac{k+1+\tau_0^{-1}}{(k+\tau_0^{-1})^2(k-1+\tau_0^{-1})}$, which is a contradiction. Therefore $\tau_k \geq \frac{1}{k+\tau_0^{-1}}$. For the other side of the inequality, assume toward contradiction that $\tau_k > \frac{2}{k+1+\tau_0^{-1}}$. Then $\frac{4(k+3+\tau_0^{-1})}{(k+1+\tau_0^{-1})^2(k-1+\tau_0^{-1})} < \frac{\tau_k^3 + \tau_k^2}{1-\tau_k} = \tau_{k-1}^2 \leq \frac{4}{(k+\tau_0^{-1})^2}$, which is a contradiction. Therefore, $\tau_k \leq \frac{2}{k+1+\tau_0^{-1}}$.

For $\{\beta_k \}$, we note that $\beta_k = \frac{\beta_{k-1}}{1 + \tau_{k-1}} = \beta_1 \prod_{i=1}^{k-1} \frac{1}{1+\tau_i} \leq \beta_1 \prod_{i=1}^{k-1} \frac{i+\tau_0^{-1}}{i+1+
\tau_0^{-1}} = \frac{\beta_1(1+\tau_0)}{\tau_0k+1}$.
\end{proof}
\end{lemma}

The following lemma is motivated by~\cite{fercoq2015accelerated}.
\begin{lemma}\label{le:g_bound}
Consider the iterates $\{\bar{x}^k, \tilde{x}^k\}_{k\geq0}$ of Algorithm~\ref{alg:A1}. Then, for $ k \geq 0$ and $i \in [n]$, we can write $\{ \bar{x}^k_i \}$ as a convex combination of $\{ \tilde{x}^l_i \}_{l=0}^k$:
\begin{equation}
\bar{x}_i^k = \sum_{l=0}^k \gamma^{k,l}_i\tilde{x}^{l}_i,
\end{equation}
where $\gamma^{k, l}_i \geq 0$ and $\sum_{l=0}^k \gamma^{k, l}_i = 1$. 
Moreover, the coefficients $\gamma^{k, l}_i$ can explicitly be computed as
\begin{align}\label{eq: gamma_coeff}
\gamma ^{k+1, l} _i = \begin{cases}
(1-\tau_k) \gamma ^{k, l} _i, & \text{for}~l = 0, \cdots, k-1, \\
(1 - \tau_k) \gamma ^{k, k}_i + \tau_k - \tfrac{\tau_k}{\tau_0}, &\text{for}~l = k, \\
\frac{\tau_k}{\tau_0}, & \text{for}~l=k+1. 
\end{cases}
\end{align}
\begin{proof}
Now, from the definition of $\bar{x}^{k+1}$ and $\hat{x}^k$, for $i\in [n]$, we can write 
\begin{equation}\label{eq:bar_x_expression}
\bar{x}^{k+1}_i = (1-\tau_k)\bar{x}^k_i +  \tau_k\tilde{x}_i^k + \frac{\tau_k}{\tau_0}(\tilde{x}^{k+1}_i - \tilde{x}^{k}_i) = (1-\tau_k)\bar{x}^k_i + (\tau_k - \frac{\tau_k}{\tau_0}) \tilde{x}_i^k +  \frac{\tau_k}{\tau_0}\tilde{x}^{k+1}_i.
\end{equation}
We prove that $\bar{x}_i^k = \sum_{l=0}^k\gamma^{k,l}_i\tilde{x}^{l}_i$ for $i \in [n]$ such that $\gamma_i^{k,l} \geq 0$ and $\sum_{l=0}^k\gamma_i^{k,l} = 1$.
Indeed, for $k=0$, we have $\bar{x}^0 = \tilde{x}^0$, which trivially holds  if we choose $\gamma_i^{0,0} = 1$. 
Now, assume that this expression holds for $k\geq 1$, we prove it holds for $k+1$.
Indeed, from \eqref{eq:bar_x_expression}, using this induction assumption, we can write
\begin{equation*}
\bar{x}_i^{k+1} = (1-\tau_k)\sum_{l=0}^{k-1}\gamma_i^{k,l}\tilde{x}_i^{l} + \bigg[(1-\tau_k)\gamma_i^{k,k} +  \tau_k - \frac{\tau_k}{\tau_0}\bigg] \tilde{x}_i^k +  \frac{\tau_k}{\tau_0} \tilde{x}^{k+1}_i = \sum_{l=0}^{k+1}\gamma^{k+1,l}_i\tilde{x}^{l}_i,
\end{equation*}
where constants $\gamma_i^{k+1,l}$ are as given in~\eqref{eq: gamma_coeff}. It is trivial to check that $\sum_{l=0}^{k+1}\gamma_i^{k+1,l} = (1-\tau_k)\sum_{l=0}^k\gamma_i^{k,l} + \tau_k - \frac{\tau_k}{\tau_0} + \frac{\tau_k}{\tau_0} = (1-\tau_k) + \tau_k = 1$.
In addition, since $\{ \tau_k \}_{k\geq0}$ is a non-increasing sequence, $\gamma_i^{k,l} \geq 0$.
\end{proof}
\end{lemma}
\section{Convergence analysis of SMART-CD}
\subsection{The proof of Theorem~\ref{th:convergence1}}
\vspace{-1ex}
First, let us define the full primal proximal-gradient step as
\begin{align}\label{eq:xtilde_bar}
\bar{\tilde{x}}^{k+1} := \argmin_{x\in\R^p}\set{\iprods{\nabla{\psi_{\beta_{k+1}}}(\hat{x}^k), x - \hat{x}^k}   + g(x) + \tau_k \sum_{i=1} ^n \frac{B^{k}_i}{2\tau_0} \Vert x_i - \tilde{x}^k_i \Vert_{(i)}^2},
\end{align}
where $\nabla \psi_{\beta_{k+1}}(\hat{x}^k) = \nabla f(\hat{x}^k) + A^\top y^\ast _{\beta_{k+1}} (A\hat{x}^k)$.
The primal coordinate step (Step~\ref{step:x_tilde}) and Step~\ref{step:x_bar} in Algorithm~\ref{alg:A1} can be written as
\begin{align}\label{eq:xtilde_bar_exp}
\tilde{x}^{k+1}_i = \begin{cases}
\bar{\tilde{x}}^{k+1}_i, & \text{if }~i = i_k,\\
\tilde{x}^{k}_i, & \text{otherwise}.
\end{cases}
\end{align}
Moreover, using \cite[Property 2]{tseng2008accelerated}, we know that for all $x \in \mathbb R^p$ and for all $i \in [n]$, 
\begin{align}\label{eq:g_bound}
g_i(\bar{\tilde{x}}^{k+1}_i) \leq  g_i(x_i) + \iprods{\nabla_i\psi_{\beta_{k+1}}(\hat{x}^k), x_i - \bar{\tilde{x}}^{k+1}_i}  &+ \frac{\tau_k B^{k}_i}{2\tau_0} \left(\Vert x_i - \tilde{x}^k_i \Vert_{(i)}^2 - \Vert x_i - \bar{\tilde{x}}^{k+1}_i \Vert_{(i)}^2 \right) \nonumber\\
&-\frac{\tau_k B^{k}_i}{2\tau_0} \Vert \bar{\tilde{x}}^{k+1}_i - \tilde{x}^k _i \Vert_{(i)}^2.
\end{align}
Now, since the partial gradient $\nabla_{i_k}{f}$ is $\hat{L}_{i_k}$-Lipschitz continuous, using $\bar{x}^{k+1}_{i_k} = \hat{x}^k_{i_k} + \frac{\tau_k}{\tau_0}(\tilde{x}^{k+1}_{i_k} - \tilde{x}^k_{i_k})$ and $\bar{x}^{k+1}_{i} = \hat{x}^k_{i}$ for $i\neq i_k$, we have
\begin{align}
f(\bar{x}^{k+1}) &\leq f(\hat{x}^k) + \iprods{\nabla _{i_{k}} f(\hat{x}^k), \bar{x}^{k+1}_{i_{k}} - \hat{x}^{k}_{i_{k}}} + \frac{\hat{L}_{i_k}}{2} \Vert \bar{x}^{k+1}_{i_{k}} - \hat{x}^{k}_{i_{k}} \Vert_{(i_k)} ^2 \nonumber\\
& = f(\hat{x}^k) + \frac{\tau_k}{\tau_0}\iprods{\nabla_{i_k}{f}(\hat{x}^k), \tilde{x}^{k+1}_{i_k} - \tilde{x}^k_{i_k}} 
 + \frac{\tau_k^2\hat{L}_{i_k}}{2\tau_0^2} \Vert \tilde{x}^{k+1}_{i_k} - \tilde{x}^{k}_{i_k} \Vert_{(i_k)}^2. 
\end{align}
Taking the $\Fc_k$-conditional expectation with respect to $i_k$ and noting~\eqref{eq:xtilde_bar_exp}, we obtain
\begin{align}
\expect{i_k}{f(\bar{x}^{k+1})\mid\Fc_k} &\leq f(\hat{x}^k) + \frac{\tau_k}{\tau_0} \sum_{i=1}^n q_i \iprods{\nabla_i f(\hat{x}^k), \bar{\tilde{x}}^{k+1}_i - \tilde{x}^{k}_i} \nonumber \\ &+ \frac{\tau_k^2}{\tau_0^2} \sum_{i=1}^n q_i\frac{\hat{L}_i}{2} \Vert \bar{\tilde{x}}^{k+1}_i - \tilde{x}^{k}_i \Vert_{(i)}^2 \label{eq: f_bound}.
\end{align}
Next, let us denote by $\varphi_{\beta}(x) := h_{\beta}(Ax)$. 
Then, by Lemma~\ref{le:smoothing_properties}, we can see that $\varphi_{\beta_{k+1}}$ has block-coordinate Lipschitz gradient with the Lipschitz constant $\frac{\Vert A_i \Vert ^2}{\beta_{k+1}}$, where $A_i$ is the $i$-th column block of $A$. 
Moreover,  $\nabla_i{\varphi_{\beta_{k+1}}}(x) = A_i^\top y^{\ast}_{\beta_{k+1}}(Ax)$.
Hence, using  $\bar{x}^{k+1}_{i_k} = \hat{x}_{i_k}^k + \frac{\tau_k}{\tau_0}(\tilde{x}^{k+1}_{i_k} - \tilde{x}^k_{i_k})$ and $\bar{x}^{k+1}_{i} = \hat{x}_{i}^k$ for $i\neq i_k$, we can write
\begin{align*}
\varphi_{\beta_{k+1}}(\bar{x}^{k+1}) &\leq \varphi_{\beta_{k+1}}(\hat{x}^k) + \iprods{\nabla_{i_k}{\varphi_{\beta_{k+1}}}(\hat{x}^k), \bar{x}_{i_k}^{k+1} - \hat{x}^k_{i_k}} + \frac{\Vert A_i\Vert^2}{2\beta_{k+1}}\Vert \bar{x}_{i_k}^{k+1} - \hat{x}^k_{i_k} \Vert_{(i_k)}^2 \nonumber\\
&= \varphi_{\beta_{k+1}}(\hat{x}^k) + \frac{\tau_k}{\tau_0}\iprods{\nabla_{i_k}{\varphi_{\beta_{k+1}}}(\hat{x}^k), \tilde{x}_{i_k}^{k+1} - \tilde{x}^k_{i_k}} + \frac{\tau_k^2\Vert A_i\Vert^2}{2\tau_0^2\beta_{k+1}}\Vert \tilde{x}_{i_k}^{k+1} - \tilde{x}^k_{i_k} \Vert_{(i_k)}^2.
\end{align*}
Taking the $\Fc_k$-conditional expectation with respect to $i_k$ given $\Fc_k$ and noting~\eqref{eq:xtilde_bar_exp}, we get
\begin{align}\label{eq:h_est1}
\expect{i_k}{\varphi_{\beta_{k+1}}(\bar{x}^{k+1})\mid \Fc_k} &\leq \varphi_{\beta_{k+1}}(\hat{x}^k)  + \frac{\tau_k}{\tau_0} \sum_{i=1}^n q_i \iprods{\nabla_i \varphi_{\beta_{k+1}}(\hat{x}^k), \bar{\tilde{x}}^{k+1}_i - \tilde{x}^k_i} \nonumber\\
& + \frac{\tau_k^2}{\tau_0^2}\sum_{i=1}^n q_i \frac{\Vert A_i\Vert^2}{2\beta_{k+1} }\Vert \bar{\tilde{x}}^{k+1}_i - \tilde{x}^k_i\Vert_{(i)}^2.
\end{align}
Now, we define 
\begin{equation}\label{eq:ghat_form}
\hat{g}_i^k := \sum_{l=0}^k\gamma_i^{k,l}g_i(\tilde{x}^l_i)~~~\text{and}~~~\hat{g}^k := \sum_{i=1}^n\hat{g}_i^k.
\end{equation}
Using Lemma~\ref{le:g_bound}, we can write 
\begin{align*}
\hat{g}^{k+1}_i &= \sum_{l=0}^{k+1}\gamma_i^{k+1,l}g_i(\tilde{x}^l_i) \\ &= \sum_{l=0}^{k-1}(1-\tau_k)\gamma_i^{k,l}g_i(\tilde{x}^l_i) + \Big[(1 - \tau_k) \gamma ^{k, k}_i + \tau_k - \tfrac{\tau_k}{\tau_0}\Big]g_i(\tilde{x}^k_i) + \frac{\tau_k}{\tau_0}g_i(\tilde{x}^{k+1}_i)\\
&= (1-\tau_k) \sum_{l=0}^k\gamma_i^{k,l}g_i(\tilde{x}^l_i) + \tau_kg_i(\tilde{x}^k_i) + \frac{\tau_k}{\tau_0}\left(g_i(\tilde{x}^{k+1}_i) - g_i(\tilde{x}_i^k)\right)\\
&= (1-\tau_k)\hat{g}_i^k + \tau_kg_i(\tilde{x}^k_i) + \frac{\tau_k}{\tau_0}\left(g_i(\tilde{x}^{k+1}_i) - g_i(\tilde{x}_i^k)\right).
\end{align*}
Using the definition \eqref{eq:ghat_form} of $\hat{g}^k$, this estimate implies
\begin{align*}
\hat{g}^{k+1} &= (1-\tau_k)\hat{g}^k + \sum_{i=1}^n\left[\tau_kg_i(\tilde{x}^k_i) + \frac{\tau_k}{\tau_0}\left(g_i(\tilde{x}^{k+1}_i) - g_i(\tilde{x}_i^k)\right)\right].
\end{align*}
Now, by the expression \eqref{eq:xtilde_bar_exp}, we can show that 
\begin{align*}
\expect{i_k}{g_i(\tilde{x}^{k+1}_i) \mid \Fc_k}  = q_ig_i(\bar{\tilde{x}}^{k+1}_i) + (1-q_i)g_i(\tilde{x}^k_i).
\end{align*}
Combining the two last expressions, we can derive
\begin{align}\label{eq:expect_of_ghat}
\expect{i_k}{\hat{g}^{k+1} \mid \Fc_k} &= (1-\tau_k)\hat{g}^k + \sum_{i=1}^n\left[\tau_kg_i(\tilde{x}^k_i) + \frac{\tau_k}{\tau_0}\left(\expect{i_k}{g_i(\tilde{x}^{k+1}_i) \mid \Fc_k} - g_i(\tilde{x}_i^k)\right)\right] \nonumber\\
&= (1-\tau_k)\hat{g}^k + \tau_k \sum_{i=1}^n g_i(\tilde{x}^k_i) + \frac{\tau_k}{\tau_0} \sum_{i=1}^nq_i \left(g_i(\bar{\tilde{x}}^{k+1}_i) - g_i(\tilde{x}_i^k)\right).
\end{align}
Let us define $\hat{F}_{\beta_k}^k := f(\bar{x}^k) + \hat{g}^k + h_{\beta_k}(A\bar{x}^k) \equiv f(\bar{x}^k) + \hat{g}^k + \varphi_{\beta_k}(\bar{x}^k)$. 
Then, from \eqref{eq: f_bound}, \eqref{eq:h_est1} and \eqref{eq:expect_of_ghat}, we have that
\begin{align}
\expect{i_k}{\hat{F}^{k+1}_{\beta_{k+1}} \mid \Fc_k} &= \expect{i_k}{f(\bar{x}^{k+1}) \mid \Fc_k} + \expect{i_k}{\hat{g}^{k+1} \mid \Fc_k} + \expect{i_k}{\varphi_{\beta_{k+1}}(\bar{x}^{k+1}) \mid \Fc_k} \nonumber\\
&\leq \left[ f(\hat{x}^k) + \frac{\tau_k}{\tau_0}\sum_{i=1}^n q_i\iprods{\nabla_i {f}(\hat{x}^k), \bar{\tilde{x}}^{k+1}_i - \tilde{x}^k_i}\right] \nonumber \\
& + \left[\varphi_{\beta_{k+1}}(\hat{x}^k) + \frac{\tau_k}{\tau_0}\sum_{i=1}^n q_i \iprods{\nabla_i \varphi_{\beta_{k+1}}(\hat{x}^k), \bar{\tilde{x}}^{k+1}_i - \tilde{x}^k_i} \right] \nonumber\\
&+ \left[(1-\tau_k)\hat{g}^k + \tau_k \sum_{i=1}^n g_i(\tilde x_i^k) +  \frac {\tau_k}{\tau_0}\sum_{i=1}^n q_i  \left(g_i(\bar{\tilde{x}}^{k+1}_i) - g_i(\tilde{x}_i^k) \right)\right] \nonumber  \\
& + \frac{\tau_k^2}{2 \tau_0^2} \sum_{i=1}^n q_i \left(\hat{L}_i + \frac{\Vert A_i\Vert^2}{\beta_{k+1}}\right) \Vert \bar{\tilde{x}}^{k+1}_i - \tilde{x}^{k}_i \Vert_{(i)}^2 \label{eq:Eik_bound},
\end{align}
since $\nabla{\psi_{\beta_{k+1}}}(\hat{x}^k) = \nabla f(\hat{x}^k) + \nabla{\varphi_{\beta_{k+1}}}(\hat{x}^k)$.
Now, using the estimate \eqref{eq:g_bound} into the last expression and noting that $B_i^k = \hat{L}_i + \frac{\Vert A_i\Vert^2}{\beta_{k+1}}$, we can further derive that for all $x$, 
\begin{align}\label{eq:before_choosing_x}
\expect{i_k}{\hat{F}^{k+1}_{\beta_{k+1}} \mid \Fc_k} 
&\leq \left[ f(\hat{x}^k) + \frac{\tau_k}{\tau_0}\sum_{i=1}^n q_i\iprods{\nabla_i {f}(\hat{x}^k), x_i - \tilde{x}^k_i}\right] \nonumber \\
& + \left[\varphi_{\beta_{k+1}}(\hat{x}^k) + \frac{\tau_k}{\tau_0}\sum_{i=1}^n q_i \iprods{\nabla_i \varphi_{\beta_{k+1}}(\hat{x}^k), x_i - \tilde{x}^k_i} \right] \nonumber\\
&+ \left[(1-\tau_k)\hat{g}^k + \tau_k \sum_{i=1}^n g_i(\tilde x_i^k) +  \frac {\tau_k}{\tau_0}\sum_{i=1}^n q_i  \left(g_i(x_i) - g_i(\tilde{x}_i^k) \right)\right] \nonumber  \\
& + \sum_{i=1}^n q_i\frac{\tau_k^2B^{k}_i}{2\tau_0^2} \big(\Vert {x}_i - \tilde{x}^k_i \Vert_{(i)}^2 - \Vert x_i - \bar{\tilde{x}}^{k+1}_i \Vert_{(i)}^2 \big).
\end{align}
Let us choose $x$ such that for all $i \in [n]$, $x_i = \left(1-\frac{\tau_0}{q_i}\right) \tilde x_i^k + \frac{\tau_0}{q_i} x_i^\star$.
Note that as $\tau_0 \leq q_i$ for all $i$, $x_i$ is a convex combination of $\tilde x_i^k$ and $x_i^{\star}$. We obtain
\begin{align}\label{eq:norms_not_simplified}
\mathbb{E}_{i_k} \Big[ \hat{F}^{k+1}_{\beta_{k+1}} \mid& \Fc_k \Big] 
\leq \left[ f(\hat{x}^k) + \tau_k \iprods{\nabla {f}(\hat{x}^k), x^\star - \tilde{x}^k}\right]  + \left[\varphi_{\beta_{k+1}}(\hat{x}^k) + \tau_k\iprods{\nabla \varphi_{\beta_{k+1}}(\hat{x}^k), x^\star - \tilde{x}^k} \right] \nonumber\\
&+ \left[(1-\tau_k)\hat{g}^k + \tau_k g(x^\star)\right] \nonumber  \\
& + \sum_{i=1}^n q_i\frac{\tau_k^2B^{k}_i}{2\tau_0^2} \left(\left\Vert\frac{\tau_0}{q_i}( x^\star_i - \tilde{x}^k_i) \right\Vert_{(i)}^2 - \left\Vert \left(1-\frac{\tau_0}{q_i}\right) \tilde x_i^k + \frac{\tau_0}{q_i} x_i^\star - \bar{\tilde{x}}^{k+1}_i \right\Vert_{(i)}^2 \right).
\end{align}
We simplify the norm difference using the fact that $\Vert ax + (1-a)y - z \Vert^2 = a \Vert x-z \Vert ^2 + (1-a) \Vert y - z \Vert ^2 - a(1-a) \Vert x-y \Vert^2$.
\begin{align*}
\bigg\Vert \left(1-\frac{\tau_0}{q_i}\right)& \tilde x_i^k + \frac{\tau_0}{q_i} x_i^\star - \bar{\tilde{x}}^{k+1}_i \bigg\Vert_{(i)}^2 \\
&=  \left(1-\frac{\tau_0}{q_i} \right) \lVert \tilde x_i^k - \bar{\tilde{x}}^{k+1}_i \rVert^2_{(i)} + \frac{\tau_0}{q_i}  \lVert x_i^\star  - \bar{\tilde{x}}^{k+1}_i \rVert^2_{(i)} - \left(1-\frac{\tau_0}{q_i}\right)\frac{\tau_0}{q_i} \lVert \tilde x_i^k - x_i^\star \rVert^2_{(i)} \\
& \geq \frac{\tau_0}{q_i}  \lVert x_i^\star  - \bar{\tilde{x}}^{k+1}_i \rVert^2_{(i)} - \left(1-\frac{\tau_0}{q_i}\right)\frac{\tau_0}{q_i} \lVert \tilde x_i^k - x_i^\star \rVert^2_{(i)}.
\end{align*}
and we get
\begin{align}\label{eq:main_est1a}
\mathbb{E}_{i_k} \Big[ \hat{F}^{k+1}_{\beta_{k+1}} \mid& \Fc_k \Big] 
\leq \big[ f(\hat{x}^k) + \tau_k \iprods{\nabla {f}(\hat{x}^k), x^\star - \tilde{x}^k}\big]  + \Big[\varphi_{\beta_{k+1}}(\hat{x}^k) + \tau_k\iprods{\nabla \varphi_{\beta_{k+1}}(\hat{x}^k), x^\star - \tilde{x}^k} \Big] \nonumber\\
&+ \left[(1-\tau_k)\hat{g}^k + \tau_k g(x^\star)\right] + \sum_{i=1}^n \frac{\tau_k^2B^{k}_i}{2\tau_0} \left(\Vert x^\star_i - \tilde{x}^k_i \Vert_{(i)}^2 - \Vert \bar{\tilde{x}}^{k+1}_i - x_i^\star \Vert_{(i)}^2 \right).
\end{align}
Using the convexity of $f$, we have $f(\hat{x}^k) + \iprods{\nabla{f}(\hat{x}^k), x^{\star} - \hat{x}^k} \leq f(x^{\star})$ and $f(\hat{x}^k) + \iprods{\nabla{f}(\hat{x}^k), \bar{x}^k - \hat{x}^k} \leq f(\bar{x}^k)$. Moreover, since $\hat{x}^k = (1-\tau_k)\bar{x}^k + \tau_k\tilde{x}^k$, we have $\tau_k(x^{\star} - \tilde{x}^k) = (1-\tau_k)(\bar{x}^k - \hat{x}^k) + \tau_k(x^{\star} - \hat{x}^k)$.
Combining these expressions, we obtain
\begin{align}\label{eq:f_est2}
f(\hat{x}^k) + \tau_k\iprods{\nabla{f}(\hat{x}^k), x^{\star} - \tilde{x}^k} &\leq (1-\tau_k)f(\bar{x}^k) + \tau_kf(x^{\star}).
\end{align}
On the one hand, by the Lipschitz gradient and convexity of $\varphi_{\beta_{k+1}}$ in Lemma~\ref{le:smoothing_properties}(b), we have 
\begin{align*}
\varphi_{\beta_{k+1}}(\hat{x}^k) + \iprods{\nabla\varphi_{\beta_{k+1}}(\hat{x}^k), \bar{x}^k - \hat{x}^k} \leq \varphi_{\beta_{k+1}}(\bar{x}^k) - \frac{\beta_{k+1}}{2}\Vert y^{\ast}_{\beta_{k+1}}(A\hat{x}^k) - y^{\ast}_{\beta_{k+1}}(A\bar{x}^k)\Vert^2.
\end{align*}
On the other hand, by Lemma~\ref{le:smoothing_properties}(c), we also have
\begin{align*}
\varphi_{\beta_{k+1}}(\hat{x}^k) + \iprods{\nabla\varphi_{\beta_{k+1}}(\hat{x}^k), x^{\star} - \hat{x}^k} \leq h(Ax^{\star}) - \frac{\beta_{k+1}}{2}\Vert y^{\ast}_{\beta_{k+1}}(A\hat{x}^k) - \dot{y}\Vert^2
\end{align*}
Combining these two inequalities and using $\tau_k(x^{\star} - \tilde{x}^k) = (1-\tau_k)(\bar{x}^k - \hat{x}^k) + \tau_k(x^{\star} - \hat{x}^k)$, we get
\begin{multline*}
\varphi_{\beta_{k+1}}(\hat{x}^k) + \tau_k\iprods{\nabla\varphi_{\beta_{k+1}}(\hat{x}^k), x^{\star} - \tilde{x}^k} \leq (1-\tau_k)\varphi_{\beta_{k+1}}(\bar{x}^k) + \tau_kh(Ax^{\star}) \\
- \frac{(1-\tau_k)\beta_{k+1}}{2}\Vert y^{\ast}_{\beta_{k+1}}(A\hat{x}^k) - y^{\ast}_{\beta_{k+1}}(A\bar{x}^k)\Vert^2 - \frac{\tau_k\beta_{k+1}}{2}\Vert y^{\ast}_{\beta_{k+1}}(A\hat{x}^k) - \dot{y}\Vert^2.
\end{multline*}
Next, using Lemma~\ref{le:smoothing_properties}(d), we can further estimate 
\begin{align}\label{eq:varphi_est2}
\varphi_{\beta_{k+1}}(\hat{x}^k) &+ \tau_k\iprods{\nabla\varphi_{\beta_{k+1}}(\hat{x}^k), x^{\star} - \tilde{x}^k} \leq (1-\tau_k)\varphi_{\beta_{k}}(\bar{x}^k) + \tau_kh(Ax^{\star}) \nonumber\\
&-\frac{(1-\tau_k)\beta_{k+1}}{2}\Vert y^{\ast}_{\beta_{k+1}}(A\hat{x}^k) - y^{\ast}_{\beta_{k+1}}(A\bar{x}^k)\Vert^2 - \frac{\tau_k\beta_{k+1}}{2}\Vert y^{\ast}_{\beta_{k+1}}(A\hat{x}^k) - \dot{y}\Vert^2 \nonumber\\
& + \frac{(1-\tau_k)(\beta_k-\beta_{k+1})}{2}\Vert y^{\ast}_{\beta_{k+1}}(A\bar{x}^k) - \dot{y}\Vert^2\nonumber\\
&\leq  (1-\tau_k)\varphi_{\beta_{k}}(\bar{x}^k) + \tau_kh(Ax^{\star}) \nonumber\\
& - \frac{1}{2}\left(\beta_{k+1}\tau_k(1-\tau_k) - (1-\tau_k)(\beta_k - \beta_{k+1})\right)\Vert y^{\ast}_{\beta_{k+1}}(A\bar{x}^k) - \dot{y}\Vert^2.
\end{align}
Here, in the last inequality, we use the fact that $(1-\tau)\Vert a - b\Vert^2 + \tau\Vert a\Vert^2 - \tau(1-\tau)\Vert b\Vert^2 = \Vert a - (1-\tau)b\Vert^2 \geq 0$ for any $a$, $b$, and $\tau\in [0, 1]$.
Substituting \eqref{eq:f_est2} and \eqref{eq:varphi_est2} into \eqref{eq:main_est1a}, we obtain
\begin{align}\label{eq:main_est2}
\expect{i_k}{\hat{F}^{k+1}_{\beta_{k+1}} \mid \Fc_k} &\leq (1-\tau_k)\left[ f(\bar{x}^k) + \hat{g}^k + \varphi_{\beta_k}(\bar{x}^k) \right] + \tau_k\left[ f(x^{\star}) + g(x^{\star}) + h(Ax^{\star})\right] \nonumber\\
&+ \sum_{i=1}^n \frac{\tau_k^2B^{k}_i}{2 \tau_0} \big(\Vert {x}^{\star}_i - \tilde{x}^k_i \Vert_{(i)}^2 - \Vert x^{\star}_i - \bar{\tilde{x}}^{k+1}_i \Vert_{(i)}^2 \big)\nonumber\\
&-\frac{(1-\tau_k)}{2}\left[\beta_{k+1}(1 + \tau_k) -  \beta_k\right]\Vert y^{\ast}_{\beta_{k+1}}(A\bar{x}^k) - \dot{y}\Vert^2.
\end{align}
Next, let us denote by $Q_k := \sum_{i=1}^n \frac{\tau_k^2B^{k}_i}{2\tau_0}\left[\Vert x^{\star}_i - \tilde{x}^{k}_i \Vert_{(i)}^2  - \Vert x^{\star}_i - \bar{\tilde{x}}^{k+1}_i \Vert_{(i)}^2 \right]$.
We can further express $Q_k$ as
\begin{align}
Q_k &= \sum_{i=1}^n\frac{\tau_k^2 B_i^k}{2\tau_0}\left[ \Vert x^{\star}_i - \tilde{x}^{k}_i \Vert_{(i)}^2  - \Vert x^{\star}_i - \bar{\tilde{x}}^{k+1}_i \Vert_{(i)}^2 \right] \nonumber \\
&= \expect{i_k}{\frac{\tau_k^2 B_{i_k}^k}{2q_{i_k} \tau_0}\left( \Vert x^{\star}_{i_k} - \tilde{x}^{k}_{i_k} \Vert_{({i_k})}^2  - \Vert x^{\star}_{i_k} - {\tilde{x}}^{k+1}_{i_k} \Vert_{({i_k})}^2 \right)\mid\Fc_k} \nonumber
\\
& = \expect{i_k}{\sum_{i=1}^n \frac{\tau_k^2 B_{i}^k}{2q_i \tau_0}\left( \Vert x^{\star}_{i} - \tilde{x}^{k}_{i} \Vert_{({i})}^2  - \Vert x^{\star}_{i} - {\tilde{x}}^{k+1}_{i} \Vert_{({i})}^2 \right)\mid\Fc_k}, \label{eq:Q_bound}
\end{align}
where the last equality follows from the fact that $\tilde{x}^{k+1}_i = \tilde{x}^k _i$ for $i \neq i_k$.

Substituting this expression into~\eqref{eq:main_est2} and using the definition of $\hat{F}_{\beta_k}^k$ and $F^\star := F(x^{\star}) = f(x^{\star}) + g(x^{\star}) + h(Ax^{\star})$, we get
\begin{align*} 
\expect{i_k}{\hat{F}^{k+1}_{\beta_{k+1}} + \sum_{i=1}^n\frac{\tau_k^2 B^{k}_{i}}{2q_{i}\tau_0}\Vert x^{\star}_{i} - \tilde{x}^{k+1}_{i} \Vert_{(i)}^2 \mid \Fc_k} &\leq (1-\tau_k)\hat{F}^k_{\beta_k} + \tau_kF(x^{\star}) \nonumber\\
&+ \sum_{i=1}^n\frac{\tau_k^2 B^{k}_{i}}{2q_{i}\tau_0} \Vert x^{\star}_{i} - \tilde{x}^{k}_{i} \Vert_{(i)}^2   - \mathcal{R}_k,
\end{align*}
where $R_k := \frac{(1-\tau_k)}{2}\left[\beta_{k+1}(1 + \tau_k) -  \beta_k\right]\Vert y^{\ast}_{\beta_{k+1}}(A\bar{x}^k) - \dot{y}\Vert^2$.
Assume that we choose $\beta_k$ and $\tau_k$ such that $\beta_{k+1}(1 + \tau_k) -  \beta_k \geq 0$, then $\mathcal{R}_k \geq 0$.
Taking the expected value of the last estimate over the $\sigma$-field $\Fc_k$, we obtain
\begin{align}\label{eq:main_est1}
\expect{}{\hat{F}^{k\!+\!1}_{\beta_{k\!+\!1}} \!-\! F^{\star}} + \expect{}{\sum_{i=1}^n\frac{\tau_k^2 B^{k}_{i}}{2q_{i}\tau_0}\Vert x^{\star}_{i} \!-\! \tilde{x}^{k\!+\!1}_{i} \Vert_{(i)}^2} &\leq (1-\tau_k)\expect{}{\hat{F}_{\beta_k}^k - F^{\star}} \nonumber\\
&+ \expect{}{\sum_{i=1}^n\frac{\tau_k^2 B^{k}_{i}}{2q_{i}\tau_0}\Vert x^{\star}_{i} - \tilde{x}^{k}_{i} \Vert_{(i)}^2}.
\end{align}
In order to telescope this inequality we assume that $\tau_k^2 B_i^k \leq (1-\tau_k)\left(\tau_{k-1}^2 B_i ^{k-1} \right)$, which is equivalent to
\begin{align}\label{eq:tau_cond}
\tau_k^2\Big(\hat{L}_i + \frac{\Vert A_i\Vert^2}{\beta_{k+1}} \Big) \leq (1-\tau_k)\bigg[\tau_{k-1}^2\Big(\hat{L}_i + \frac{\Vert A_i\Vert^2}{\beta_{k}} \Big) \bigg].
\end{align}
Let us update $\beta_{k+1} = \frac{\beta_k}{1 + \tau_k}$. 
Then, this condition becomes
\begin{align}\label{eq:param_condition} 
\tau_k^2\left(\beta_k\hat{L}_i \!+\! (1 \!+\! \tau_k)\Vert A_i\Vert^2 \right) \leq (1 \!-\! \tau_k)\tau_{k-1}^2\left(\beta_k\hat{L}_i + \Vert A_i\Vert^2 \right).
\end{align}
The condition~\eqref{eq:param_condition} holds if $\tau_k^2(1+\tau_k) = (1-\tau_k)\tau_{k-1}^2$. Hence, we can compute $\tau_k$ as the unique positive root of $\tau^3 + \tau^2 + \tau_{k-1}^2\tau - \tau_{k-1}^2 = 0$.
By Lemma~\ref{le:parameters}, the root of this cubic satisfies $\frac{1}{k+\tau_0^{-1}} \leq \tau_k \leq \frac{2}{k+\tau_0^{-1} + 1}$.
Let us define $S_k = \sum_{i=1}^n \frac{\tau_k^2 B_i ^k}{2q_i \tau_0} \Vert x_i ^\star - \tilde{x}^{k+1}_i \Vert _{(i)}^2$. Then, we can recursively show that 
\begin{align*}
\expect{}{\hat{F}^{k+1}_{\beta_{k+1}} - F^{\star} + S_k} &\leq \prod_{i=1}^k(1-\tau_i)\expect{}{\hat{F}_{\beta_1}^1 - F^{\star} + \sum_{i=1}^n\frac{\tau_0^2B^{0}_{i}}{2q_{i}\tau_0} \Vert x^{\star}_{i} - \tilde{x}^{1}_{i} \Vert_{(i)}^2} \\
& \leq \prod_{i=1}^k(1-\tau_i)\left((1-\tau_0)(\hat{F}_{\beta_0}^0 - F^{\star}) + \sum_{i=1}^n\frac{\tau_0^2B^{0}_{i}}{2q_{i}\tau_0} \Vert x^{\star}_{i} - \tilde{x}^{0}_{i} \Vert_{(i)}^2\right),
\end{align*}
where the second inequality follows from~\eqref{eq:main_est1}. 
Since $\tau_k \geq \frac{1}{k + \tau_0^{-1}}$, it is trivial to show that $\omega_{k+1} := \prod_{i=1}^k (1-\tau_i)\leq \prod_{i=1}^k\frac{i + \tau_0^{-1}-1}{i+\tau_0^{-1}} = \frac{1}{\tau_0 k + 1}$.
Now, we have $F_{\beta_0}(x^0) = f(x^0) + g(x^0) + h_{\beta_0}(Ax^0) = \hat{F}_{\beta_0}^0$, and $\tilde{x}^0 = x^0$.  
In addition, by the convexity of $g$ and Lemma~\ref{le:g_bound}, we also have $g(\bar{x}^k) = g\left(\sum_{l=0}^k\gamma^{k,l}\tilde{x}^{l}\right) \leq \sum_{l=0}^k\gamma^{k,l}g(\tilde{x}^l) = \hat{g}^k$.
Hence, we can write the above estimate as
\begin{equation}
\label{eq:smoothed_gap_bound}
\expect{}{F_{\beta_k}(\bar{x}^k) - F^{\star}} \leq \frac{1}{\tau_0(k-1) + 1}\left[(1-\tau_0)(F_{\beta_0}(x^0) - F^{\star}) + \sum_{i=1}^n\frac{\tau_0B^{0}_{i}}{2q_{i}} \Vert x^{\star}_{i} - {x}^{0}_{i} \Vert_{(i)}^2\right].
\end{equation}
Now, using the bound \eqref{eq:bound_est}, we have $0 \leq F(\bar{x}^k) - F^{\star} \leq F_{\beta_k}(\bar{x}^k) - F^{\star} + \beta_k\frac{D_{h^{\ast}}^2}{2}$. Combining this estimate and the above inequality, and noting that $\beta_k \leq \frac{\beta_1(1+\tau_0)}{\tau_0k+1}$ by Lemma~\ref{le:parameters}, we obtain the bound in \eqref{eq:main_result1}.
\Eproof

\subsection{The proof of Theorem~\ref{th:convergence2}}
\vspace{-2ex}
Since $h(u) = \delta_{\set{c}}(u)$, we can smooth this function as $h_{\beta}(u) = \max_y\set{\iprods{u-c, y} - \frac{\beta}{2}\Vert y - \dot{y}\Vert^2}$.
Let us first define $S_{\beta}(x) := \expect{}{F(x) + h_{\beta}(Ax) - F^{\star}}$. 
Since $h^{\ast}(y) = \iprods{c, y}$, we use Lemma~\ref{le:smoothing_properties}(e) to estimate \eqref{eq:varphi_est2} in the proof of Theorem~\ref{th:convergence1} instead of Lemma~\ref{le:smoothing_properties}(d) to obtain  
\begin{align}\label{eq:varphi_est2_b}
\varphi_{\beta_{k+1}}(\hat{x}^k) &+ \tau_k\iprods{\nabla\varphi_{\beta_{k+1}}(\hat{x}^k), x^{\star} - \tilde{x}^k} \leq  (1-\tau_k)\varphi_{\beta_{k}}(\bar{x}^k) + \tau_kh(Ax^{\star}) \nonumber\\
& - \frac{(1-\tau_k)\beta_{k+1}}{2\beta_k}\left[\beta_{k+1} - (1-\tau_k)\beta_k\right]\Vert y^{\ast}_{\beta_{k+1}}(A\bar{x}^k) - \dot{y}\Vert^2.
\end{align}
Hence, if $\beta_{k+1} = (1-\tau_k)\beta_k$, then $\varphi_{\beta_{k+1}}(\hat{x}^k) + \tau_k\iprods{\nabla\varphi_{\beta_{k+1}}(\hat{x}^k), x^{\star} - \tilde{x}^k} \leq  (1-\tau_k)\varphi_{\beta_{k}}(\bar{x}^k) + \tau_kh(Ax^{\star})$.
Now, we combine the condition $\beta_{k+1} = (1-\tau_k)\beta_k$ and \eqref{eq:tau_cond}, we can show that 
\begin{equation*}
\tau_k^2\left((1-\tau_k)\beta_k\hat{L}_i + \Vert A_i\Vert^2\right) \leq (1-\tau_k)^2\tau_{k-1}^2\left(\beta_k\hat{L}_i + \Vert A_i\Vert^2\right).
\end{equation*}
This condition holds if $\tau_k^2 = (1-\tau_k)^2\tau_{k-1}^2$, which leads to $\tau_k = \frac{\tau_{k-1}}{\tau_{k-1} + 1}$. This is the update rule \eqref{eq:update_tau_beta} of the algorithm. It is trivial to show that $\tau_k = \frac{1}{k + \tau_0^{-1}}$ and $\beta_k = \frac{\beta_1}{\tau_0(k-1) + 1}$.
Now, we apply \eqref{eq:smoothed_gap_bound} to obtain the bound
\begin{equation*}
S_{\beta_k}(\bar{x}^k) \leq \frac{C^{\ast}}{\tau_0(k-1)+1},~~~\text{where}~~C^{\ast} :=  (1-\tau_0)(F_{\beta_0}(x^0) - F^{\star}) + \sum_{i=1}^n{\!\!}\frac{\tau_0B^{0}_{i}}{2q_{i}} \Vert x^{\star}_{i} - {x}^{0}_{i} \Vert_{(i)}^2.
\end{equation*}
Now, let us define the dual problem of~\eqref{eq:cvx_prob} as
\begin{equation}\label{eq:dual_prob}
	\max_{y \in \mathbb{R}^m} \left\{ \min_{x\in\mathbb{R}^p} F(x) + \langle Ax, y \rangle - h^\ast(y) \right\},
\end{equation}
and denote an optimal point of~\eqref{eq:dual_prob} as $y^\star$. We define $D_{\beta_k}(x) := F(x) + h_{\beta_k}(Ax) - F^\star$ and apply \cite[Lemma~1]{tran2015smooth} to obtain algorithm-independent duality bounds
\begin{align}
\left\{\begin{array}{lll}\label{eq:const_gen_bounds}
&F(\bar{x}^k) - F^{\star} &\leq D_{\beta_k}(\bar x^k) + \norm{y^\star} \norm{A \bar x^k - b} + \frac{\beta_k}{2} \norm{y^\star - \dot y}^2, \vspace{1.5ex} \\
&\Vert A\bar{x}^k - b\Vert &\leq \beta_k \left[\Vert y^{\star}-\dot{y}\Vert + \left(\Vert y^{\star}-\dot{y}\Vert^2 + 2\beta_k^{-1} D_{\beta_k}(\bar x^k)\right)^{1/2}\right].\end{array}\right.
\end{align}
The result in~\eqref{eq:main_result2} follows by taking the expectation and using the concavity of the square-root and Jensen's inequality.
\Eproof

\subsection{The proof of Corollary~\ref{co:zero_g_case}}
\vspace{-1ex}
From the update in~\eqref{eq:coordinate_grad_step}, we get the trivial inequality, similar to~\eqref{eq:g_bound}, that
\begin{align}
\frac{\tau_k B_i^k}{2q_i} \Vert \bar{\tilde{x}}^{k+1}_i - \tilde{x}^k_i \Vert_{(i)}^2 &\leq \langle \nabla _{i} \psi_{\beta_{k+1}}(\hat{x}^k), x^\star_i - \bar{\tilde{x}}^{k+1}_i \rangle \nonumber \\ &+ \frac{\tau_k B_i^k}{2q_i} \left( \Vert \bar{\tilde{x}}^{k+1}_i - x^\star_i \Vert_{(i)}^2 -\Vert{\tilde{x}}^{k}_i - x^\star_i \Vert_{(i)}^2 \right) \label{eq:g0_first}.
\end{align}
Due to the specific Step~\ref{step:x_bar} in Section~\ref{subsec:spec_cases}, instead of~\eqref{eq:Eik_bound}, we get
\begin{align}
\expect{i_k}{\hat{F}^{k+1}_{\beta_{k+1}} \mid \Fc_k} &= \expect{i_k}{f(\bar{x}^{k+1}) \mid \Fc_k} + \expect{i_k}{\hat{g}^{k+1} \mid \Fc_k} + \expect{i_k}{\varphi_{\beta_{k+1}}(\bar{x}^{k+1}) \mid \Fc_k} \nonumber\\
&\leq \left[ f(\hat{x}^k) + {\tau_k}\sum_{i=1}^n \iprods{\nabla_i {f}(\hat{x}^k), \bar{\tilde{x}}^{k+1}_i - \tilde{x}^k_i}\right] \nonumber \\
& + \left[\varphi_{\beta_{k+1}}(\hat{x}^k) + \tau_k\sum_{i=1}^n \iprods{\nabla_i \varphi_{\beta_{k+1}}(\hat{x}^k), \bar{\tilde{x}}^{k+1}_i - \tilde{x}^k_i} \right] \nonumber\\
& + \sum_{i=1}^n \frac{\tau_k^2}{2q_i} \left(\hat{L}_i + \frac{\Vert A_i\Vert^2}{\beta_{k+1}}\right) \Vert \bar{\tilde{x}}^{k+1}_i - \tilde{x}^{k}_i \Vert_{(i)}^2.
\end{align}
Plugging\eqref{eq:g0_first} into the last inequality gives us
\begin{align}
\expect{i_k}{\hat{F}^{k+1}_{\beta_{k+1}} \mid \Fc_k} &= \expect{i_k}{f(\bar{x}^{k+1}) \mid \Fc_k} + \expect{i_k}{\hat{g}^{k+1} \mid \Fc_k} + \expect{i_k}{\varphi_{\beta_{k+1}}(\bar{x}^{k+1}) \mid \Fc_k} \nonumber\\
&\leq \left[ f(\hat{x}^k) + {\tau_k} \langle \nabla f(\hat{x}^k), x^\star - \tilde{x}^k \rangle \right] + \left[\varphi_{\beta_{k+1}}(\hat{x}^k) + \tau_k \langle \nabla\varphi_{\beta_{k+1}}(\hat{x}^k), x^\star - \tilde{x}^k \rangle \right] \nonumber\\
& + \sum_{i=1}^n \frac{\tau_k^2 B_i^k}{2q_i} \left( \Vert x^\star _i - \tilde{x}^{k}_i \Vert _{(i)}^2 - \Vert x^\star_i - \bar{\tilde{x}}^{k+1}_i \Vert_{(i)}^2 \right).
\end{align}
If we let $Q_k := \sum_{i=1}^n \frac{\tau_k^2 B_i^k}{2q_i} \left( \Vert x^\star _i - \tilde{x}^{k}_i \Vert _{(i)}^2 - \Vert x^\star_i - \bar{\tilde{x}}^{k+1}_i \Vert_{(i)}^2 \right)$, then similar to~\eqref{eq:Q_bound}, we get
\begin{align}
Q_k = \expect{i_k}{\sum_{i=1}^n \frac{\tau_k^2 B_i^k}{2q_i^2} \left( \Vert x^\star_i - \tilde{x}^k_i \Vert ^2 _{(i)} - \Vert x^\star_i - \tilde{x}^{k+1}_i \Vert ^2 _{(i)} \right)\mid\Fc_k}.
\end{align}
Consequently, by using the same updates for $\tau_k$ and $\beta_k$, the recursion in~\eqref{eq:main_est1} becomes
\begin{align}
\expect{}{\hat{F}^{k\!+\!1}_{\beta_{k\!+\!1}} \!-\! F^{\star}} + \expect{}{\sum_{i=1}^n\frac{\tau_k^2 B^{k}_{i}}{2q_{i}^2}\Vert x^{\star}_{i} \!-\! \tilde{x}^{k\!+\!1}_{i} \Vert_{(i)}^2} &\leq (1-\tau_k)\expect{}{\hat{F}_{\beta_k}^k - F^{\star}} \nonumber\\
&+ \expect{}{\sum_{i=1}^n\frac{\tau_k^2 B^{k}_{i}}{2q_{i}^2}\Vert x^{\star}_{i} - \tilde{x}^{k}_{i} \Vert_{(i)}^2}.
\end{align}
Hence, we finally get
\begin{equation*}
\expect{}{F_{\beta_k}(\bar{x}^k) - F^{\star}} \leq \frac{1}{\tau_0(k-1) + 1}\left[(1-\tau_0)(F_{\beta_0}(x^0) - F^{\star}) + \sum_{i=1}^n\frac{\tau_0^2 B^{0}_{i}}{2q_{i}^2} \Vert x^{\star}_{i} - {x}^{0}_{i} \Vert_{(i)}^2\right].
\end{equation*}
Noting that $\tau_0=1$, we get
\begin{equation}\label{eq:S_beta_spec}
S_{\beta_k}(\bar{x}^k) \leq \frac{C^{\ast}}{k},~~~\text{where}~~C^{\ast} := \sum_{i=1}^n{\!\!}\frac{B^{0}_{i}}{2q_{i}^2} \Vert x^{\star}_{i} - {x}^{0}_{i} \Vert_{(i)}^2.
\end{equation}
By using the bound of $\{ \beta_k \}_{k\geq1}$ as in Lemma~\ref{le:parameters}, we obtain the bound~\eqref{eq:main_result1c}.
For the constrained case, we use~\eqref{eq:S_beta_spec} on~\eqref{eq:const_gen_bounds}, with the specific update rule of $\{ \beta_k \}$ for the constrained case, to obtain~\eqref{eq:main_result2c} using the same arguments  as in the Proof of Theorem~\ref{th:convergence2}.
\Eproof
\section{Equivalence of SMART-CD and Efficient SMART-CD}
In this appendix, we give a proof by induction for the equivalence of Algorithm~\ref{alg:A1} and Algorithm~\ref{alg:A_eff} motivated by~\cite{fercoq2015accelerated}.
\subsection{The proof of Proposition~\ref{prop:eff_equiv}}
\vspace{-1ex}
The claim trivially holds for $k=0$ using the initialization of the parameters. Assume that the relations hold for some $k$. Using Step ~\ref{step:z_tilde} of Algorithm~\ref{alg:A_eff}, we have
\begin{equation}\label{eq:zt_update}
\tilde{z}^{k+1}_{i_k} = \tilde{z}^{k}_{i_k} + t^{k+1}_{i_k}.
\end{equation}
We can write from Step~\ref{step:t} of Algorithm~\ref{alg:A_eff} that
\begin{align*}
t^{k+1} _{i_k} &= \mathrm{arg}\!\!\!\min_{t\in\R^{p_{i_k}}} \bigg\{ \langle \nabla _{i_k} f(c_k u^k + \tilde{z}^k) + A^{\top}_{i_k} y_{\beta_{k+1}}^{\ast} \big(c_k Au^k + A\tilde{z}^k\big), t \rangle + g_{i_k}(t + \tilde{z}^k_{i_k}) \nonumber \\ &~~~~~~~~~~~~~~~~~~~~~~+ \frac{\tau_k B^k _{i_k}}{2\tau_0} \Vert t \Vert_{(i_k)} ^2 \bigg\} \nonumber\\
&= \mathrm{arg}\!\!\!\min_{t\in\R^{p_{i_k}}} \bigg\{ \langle \nabla _{i_k} f (\hat{z}^k) + A^{\top}_{i_k} y_{\beta_{k+1}}^{\ast} \big(A\hat{z}^k\big), t \rangle + g_{i_k}(t + \tilde{z}^k_{i_k}) + \frac{\tau_k B^k _{i_k}}{2\tau_0} \Vert t \Vert ^2 _{(i_k)} \bigg\} \nonumber\\
&= \mathrm{arg}\!\!\!\min_{t\in\R^{p_{i_k}}} \bigg\{ \langle \nabla _{i_k} f (\hat{x}^k) + A^{\top}_{i_k} y_{\beta_{k+1}}^{\ast} \big(A\hat{x}^k\big), t \rangle + g_{i_k}(t + \tilde{x}^k_{i_k}) + \frac{\tau_k B^k _{i_k}}{2\tau_0} \Vert t \Vert ^2 _{(i_k)} \bigg\} \nonumber\\
&= - \tilde{x}^{k} _{i_k} + \mathrm{arg}\!\!\!\min_{x\in\R^{p_{i_k}}} \bigg\{ \langle \nabla _{i_k} f (\hat{x}^k) + A^{\top}_{i_k} y_{\beta_{k+1}}^{\ast} \big(A\hat{x}^k\big), x - \hat{x}^k _{i_k} \rangle + g_{i_k}(x) \nonumber \\ &~~~~~~~~~~~~~~~~~~~~~~~~~~~~~~~~~~~~~+ \frac{\tau_k B^k _{i_k}}{2\tau_0} \Vert x - \tilde{x}^k_{i_k} \Vert ^2 _{(i_k)}\bigg\} \nonumber\\
&= - \tilde{x}^{k} _{i_k} + \tilde{x}^{k+1}_{i_k}.
\end{align*}
By \eqref{eq:zt_update} and the inductive assumption on $\tilde{x}^k$, we obtain
\begin{equation*}\label{eq:equiv_tilde}
\tilde{z}^{k+1} = \tilde{x}^{k+1}.
\end{equation*}
Next, using the definition of $\bar{z}^{k+1}$ and Step~\ref{step:u}, we can derive
\begin{align*}
\bar{z}^{k+1} &= c_k u^{k+1} + \tilde{z}^{k+1} = c_k \bigg(u_k - \frac{1 - \tau_k/\tau_0}{c_k} (\tilde{z}^{k+1} - \tilde{z}^k) \bigg) + \tilde{z}^{k+1} \nonumber\\
&= c_k u^k + \tilde{z}^k + \frac{\tau_k}{\tau_0} (\tilde{z}^{k+1} - \tilde{z}^k) \nonumber\\
&= \hat{z}^k + \frac{\tau_k}{\tau_0} (\tilde{z}^{k+1} - \tilde{z}^k) \nonumber\\
&= \hat{x}^k + \frac{\tau_k}{\tau_0} (\tilde{x}^{k+1} - \tilde{x}^k) \nonumber\\
&= \bar{x}^{k+1}. 
\end{align*}
Finally, we use the definition of $\hat{z}^{k+1}$, $c_k$ and Step~\ref{step:update_tau_eff} of Algorithm~\ref{alg:A1}, we arrive at
\begin{align*}
\hat{z}^{k+1} &= c_{k+1} u^{k+1} + \tilde{z}^{k+1} \\
&= \frac{c_{k+1}}{c_k} (\bar{x}^{k+1} - \tilde{z}^{k+1}) + \tilde{z}^{k+1} \\
&= (1 - \tau_{k+1}) (\bar{z}^{k+1} - \tilde{z}^{k+1}) + \tilde{z}^{k+1} \\
&= (1 - \tau_{k+1}) (\bar{x}^{k+1} - \tilde{x}^{k+1}) + \tilde{x}^{k+1} \\
&= (1 - \tau_{k+1}) \bar{x}^{k+1} + \tau_{k+1} \tilde{x}^{k+1} \\
&= \hat{x}^{k+1}.
\end{align*} 
Hence, we can conclude that Algorithm~\ref{alg:A1} and Algorithm~\ref{alg:A_eff} are equivalent.
\Eproof

\end{document}